\documentclass[11pt,a4paper]{article}
\usepackage{amssymb}
\usepackage{amsmath}
\usepackage{latexsym}
\usepackage{graphicx} 
\usepackage{hyperref}
\usepackage{amsthm}
\usepackage{verbatim}
\usepackage{psfrag}
\usepackage{bbm}
\usepackage{dsfont}
\usepackage{mathrsfs}
\usepackage{tikz,pgfplots}
\usepackage{yfonts}
\usepackage{color}
\usepackage{enumerate}
\usepackage{subcaption}
\usepackage{caption}
\usepackage{pgfplots}
\usepackage{subcaption}
\usepackage{multirow}
\usepackage{titling}

\newcommand*{\collauthor}[2]{{#1}$^{#2}$}

\newcommand*{\affiliation}[2]{$\mbox{}^{{#2}}${#1}}
\newcommand*{\colltitle}[1]{\textbf{#1}}

\newtheorem{mydef}{Definition}[section]
\newtheorem{pro}{Proposition}[section]
\newtheorem{lemma}{Lemma}[section]
\newtheorem{theorem}[lemma]{Theorem}
\newtheorem{scheme}[]{Scheme}
\newtheorem{remark}[pro]{Remark}
\newenvironment{keywords}[1]{\vspace{1cm}\\{\bf \slshape{Keywords}}\quad\slshape{#1}}{}

\newcommand{\E}{\mathds{E}}

\newcommand*{\tbi}[3]{{#1}_{#2}^{#3}}

\newenvironment{alphafootnotes}
  {\par\edef\savedfootnotenumber{\number\value{footnote}}
   
   \setcounter{footnote}{0}}
  {\par\setcounter{footnote}{\savedfootnotenumber}}

\begin{document}
\begin{center}
\begin{Large}
  \colltitle{High-order combined Multi-step Scheme for solving forward Backward Stochastic Differential Equations}
\end{Large} 
\vspace*{1.5ex}
\begin{alphafootnotes}
\begin{sc}
\begin{large}
\collauthor{Long Teng}{1,}\footnote{Corresponding author (teng@math.uni-wuppertal.de)}
and
\collauthor{Weidong Zhao}{2}
\end{large}
\end{sc}
\end{alphafootnotes}
\vspace{1.5ex}

\affiliation{Lehrstuhl f\"ur Angewandte Mathematik und Numerische Analysis,\\
Fakult\"at f\"ur Mathematik und Naturwissenschaften,\\
Bergische Universit\"at Wuppertal, Gau{\ss}str. 20, 
42119 Wuppertal, Germany}{1} \\
\affiliation{School of Mathematics \& Finance Institute,\\ Shandong University, Jinan 250100, China(wdzhao@sdu.edu.cn)}{2} 
\end{center}

\section*{Abstract}
In this work, in order to obtain higher-order schemes for solving forward backward stochastic differential equations, we adopt the high-order multi-step method in [W. Zhao, Y. Fu and T. Zhou, SIAM J. Sci. Comput., 36(4) (2014), pp.A1731-A1751] by combining multi-steps. Two reference ordinary differential equations containing the conditional expectations and their derivatives are derived from the backward component. These derivatives are approximated by finite difference methods with multi-step combinations. The resulting scheme is a semi-discretization in the time direction involving conditional expectations, which are solved by using the Gaussian quadrature rules and polynomial interpolations on the spatial grids. Our new proposed multi-step scheme allows for higher convergence rate up to ninth order, and are more efficient. Finally, we provide a numerical illustration of the convergence of the proposed method.
\begin{keywords}
forward backward stochastic differential equations, multi-step scheme, finite difference method, time-space grid, Gauss-Hermite quadrature rule
\end{keywords}

\section{Introduction}
Recently, the forward-backward stochastic differential equation (FBSDE) becomes an important tool for formulating many problems in various areas including physics and financial mathematics. We are interested in the numerical approximation of the general FBSDEs
\begin{equation}\label{eq:fbsde}
\left\{
\begin{array}{l}
\,\,\, dX_t = a(t, X_t, Y_t, Z_t)\,dt + b(t, X_t, Y_t, Z_t)\,dW_t,\quad X_0=x_0, \hspace*{0.69cm} \mbox{forward component}\\
-dY_t = f(t, X_t, Y_t, Z_t)\,dt - Z_t\,dW_t,\hspace*{4.49cm}\mbox{backward component}\\
\quad Y_T=\xi=g(X_T)
 \end{array}\right.
 \end{equation}
on a filtered complete probability space $(\Omega, \mathcal{F}, P)$ with the natural filtration $(\mathcal{F}_t)_{0\leq t \leq T},$ where $a:[0, T] \times \mathbb{R}^n \times \mathbb{R}^m \times \mathbb{R}^{m\times d}  \to \mathbb{R}^n$ and $b:[0, T] \times \mathbb{R}^n \times \mathbb{R}^m \times \mathbb{R}^{m\times d}  \to \mathbb{R}^{n\times d},$ are drift and diffusion coefficients in the forward component, respectively; $W_t=(W^1_t,\cdots, W^d_t)^T$ is a $d$-dimensional Brownian motion (all Brownian motions are independent with each other); $f(t, X_t, Y_t, Z_t): [0, T] \times \mathbb{R}^n \times \mathbb{R}^m \times \mathbb{R}^{m\times d} \to \mathbb{R}^m$ is the driver function and $\xi$ is the square-integrable terminal condition. We see that the terminal condition $Y_T$ depends on final value of the forward component. Note that $a, b$ and $f$ are all $\mathcal{F}_t$-adapted, and a triple $(X_t, Y_t, Z_t)$ is called an $L^2$-adapted solution of \eqref{eq:fbsde} if it is $\mathcal{F}_t$-adapted, square integrable, and satisfies
\begin{equation}
\left\{
\begin{array}{l}
X_t = X_0 + \int_0^t a(s, X_s, Y_s, Z_s)\,ds + \int_0^t b(s, X_s, Y_s, Z_s)\,dW_s,\\
Y_t = \xi + \int_t^T f(s, X_s, Y_s, Z_s)\,ds - \int_t^T Z_s\,dW_s.
\end{array}\right.
\end{equation} 
One obtains decoupled FBSDEs if $a$ and $b$ are independent with $Y_t$ and $Z_t$ in \eqref{eq:fbsde}, which become
backward stochastic differential equations (BSDEs) when $a=0$ and $b=1.$

The existence and uniqueness of solution of the BSDEs assuming the Lipschitz conditions on $f,a,b~\mbox{and}~g$ are proven by Pardoux and Peng \cite{Pardoux1990, Pardoux1992}. The uniqueness of solution is extended under more general assumptions for $f$ in \cite{Lepeltier1997},
but only in the one-dimensional case. The existence and uniqueness of solution of FBSDEs have been studied in \cite{Ma1994, Peng1999}.

In recent years, many numerical methods have been proposed for the BSDEs and FBSDEs. We list some of them here: \cite{Bender2008, Bender2012, Bouchard2004, Crisan2010, Cvitanic2006, Delarue2006, Douglas1996, Fu2017, Gobet2005, Lemor2006, Ma2005, Ma2008, Milstein2006, Ruijter2015, Teng2019, Teng2020, Zhao2006, Zhao2009, Zhao2010, Zhao2012, Zhao2013, Zhao2014a, Zhang2004, Zhang2013}, and many others. In this literature, the high-order methods rely on the high-order approaches for both the forward and backward components, where are clearly difficult and computationally expensive to achieve.

Moreover, Zhao et al. proposed in \cite{Zhao2014} new kinds of high-order multi-step schemes for FBSDEs, which can keep high-order accuracy while using the Euler method to solve the forward component. This is of great interesting since the use of Euler method can dramatically simplify the entire computations. However,
the convergence rate is restricted to sixth order, since the stability condition cannot be satisfied for a higher order. For this reason, we adopt in this work this method by combining some multi-steps to achieve higher rate of convergence. More precisely, we derive two reference ordinary differential equations (ODEs) which contain the conditional expectations and their derivatives. We approximate these derivatives using finite difference methods with the combination of multi-steps for a better stability. The resulting conditional expectations are solved using the Gaussian quadrature rules, whereas the first component, $X_t$ is solved using the Euler method thanks to the local property of the generator of diffusion processes. FBSDEs are numerically solved on the time-space grids.
Numerical experiments are presented to demonstrate the improvement in the rate of convergence.

In the next section, we start with preliminaries on FBSDEs and derive in Section 3 the approximations of derivatives using finite difference method with combined multi-steps. In Section 4, we derive the reference ODEs, based on which the semi-discrete higher-order multi-step schemes are introduced for solving decoupled FBSDEs. Section 5 is devoted to the fully discrete higher-order schemes. In Section 6, these methods are extended to solve a coupled FBSDE. In Section 7, several numerical experiments on the decoupled and coupled FBSDEs including two-dimensional applications are provided to show the higher efficiency and accuracy. Finally, Section 8 concludes this work.

\section{Preliminaries}\label{sec:Preliminaries}
As mentioned before, throughout the paper we assume that $(\Omega, \mathcal{F}, P)$ is a complete, filtered probability space. A standard $d$-dimensional Brownian motion $W_t$ with a finite terminal time $T$ is defined, and the first component, $X_t$ generates the filtration $\mathcal{F}_t=\sigma\{X_s, 0\leq s \leq t\}.$ And the usual hypotheses should be satisfied. We denote the set of all $\mathcal{F}_t$-adapted and square
integrable processes in $\mathbb{R}^d$ with $L^2=L^2(0,T;\mathbb{R}^d),$ and list following notation to be used:
\begin{itemize}
	\item $|\cdot|:$ the Euclidean norm in $\mathbb{R},$ $\mathbb{R}^n$ and $\mathbb{R}^{n\times d};$
	\item $\mathcal{F}_t^{s,x}:$ $\sigma$-algebra generated by the diffusion process $\{X_r, s\leq r\leq t, X_s=x\};$
	\item $\tbi{\E}{t}{s,x}[\cdot]:$ conditional expectation under $\mathcal{F}_t^{s,x},$ i.e., $\tbi{\E}{t}{s,x}[\cdot|\mathcal{F}_t^{s,x}];$
	\item $C_b^k:$ the set of continuous functions with uniformly bounded derivatives up to order $k;$
	\item $C^{k_1,k_2}:$ the set of functions with continuous partial derivatives $\frac{\partial}{\partial t}$ and $\frac{\partial}{\partial x}$ up to $k_1$ and $k_2,$ respectively;
	\item $C_L:$ the set of uniformly Lipschitz continuous function with respect to the spatial variables;
	\item $C^{\frac{1}{2}}_L:$ the subset of $C_L$ such that its element is H\"older-$\frac{1}{2}$ continuous with respect to time, with  uniformly bounded Lipschitz and H\"older constants.
\end{itemize}
Let $X_t$ be a diffusion process 
\begin{equation}\label{eq:DiffusionX}
X_t = x_0 + \int_0^t a(s, X_s)\,ds + \int_0^t b(s, X_s)\,dW_s
\end{equation}
starting at $(t_0, x_0)$ and $t\in[t_0, T],$ which has a unique solution. Note that  $\tbi{\E}{s}{x}[X_t]:=\tbi{\E}{s}{s,x}[X_t]$ is equal to
$\E[X_t|X_s=x]$ for all $s \leq t$ with the Markov property of the diffusion process. Given a measurable function $g:[0, T] \times  \mathbb{R}^n \to \mathbb{R},$ $\tbi{\E}{s}{x}[g(t, X_t)]$ is a function of $(t,s,x)$ whose partial derivative with respect to $t$ reads
\begin{equation*}
	\frac{\partial \tbi{\E}{s}{x}[g(t, X_t)]}{\partial t}=\lim_{\tau \to 0^+}\frac{\tbi{\E}{s}{x}[g(t+\tau, X_{t+\tau})]-\tbi{\E}{s}{x}[g(t, X_t)]}{\tau}
\end{equation*}
provided that the limit exists and is finite.
\begin{mydef}[Generator]The generator $\mathcal{A}_t^x$ of $X_t$ satisfying \eqref{eq:DiffusionX} on a measurable function $g:[0, T] \times  \mathbb{R}^n \to \mathbb{R}$ is defined by
\begin{equation*}
\mathcal{A}_t^xg(t,x)=\lim_{h \to 0^+}\frac{\tbi{\E}{t}{x}[g(t+h, X_{t+h})]-g(t, x)}{h},\quad x\in\mathbb{R}^n.
\end{equation*}	
\end{mydef}
\begin{theorem}\label{theorem:Generator}
	Let $X_t$ be the diffusion process defined by \eqref{eq:DiffusionX}, then it holds
	  \begin{equation}\label{eq:GeneratorX}
	  	\mathcal{A}_t^x f(t, x)=\mathcal{L}_{t, x}f(t, x)
	  \end{equation}
	  for $f \in C^{1,2}\left([0,T]\times \mathbb{R}^n\right)$ with
	  \begin{equation*}
	  	\mathcal{L}_{t,x}=\frac{\partial}{\partial t}+\displaystyle\sum_{i}a_i(t,x)\frac{\partial}{\partial x_i}+\frac{1}{2}\displaystyle\sum_{i,j}(bb^{\top})_{i,j}(t,x)\frac{\partial^2}{\partial x_i \partial x_j}.
	  \end{equation*}
\end{theorem}
The proof can be simply completed by using the It\^o's lemma and the dominated convergence theorem.
\begin{remark}
	From \eqref{eq:GeneratorX} one can straightforwardly deduce that
	 \begin{equation*}
	\mathcal{A}_t^{X_t} f(t, X_t)=\mathcal{L}_{t, X_t}f(t, X_t),
	\end{equation*}
	which is a stochastic process.
\end{remark}
By using the It\^o's lemma and Theorem \ref{theorem:Generator} we calculate
\begin{equation*}
		\left.\frac{d\,\tbi{\E}{t_0}{x_0}[f(t, X_t)]}{d t}\right|_{t=t_0}=\lim_{t \to t_0^+}\frac{\tbi{\E}{t_0}{x_0}[f(t, X_{t})]-g(t_0, x_0)}{t-t_0}=\mathcal{L}_{t, x}f(t_0, x_0)=\mathcal{A}_t^{x} f(t_0, x_0),
\end{equation*}
from which we deduce Theorem \eqref{theorem:derivatives} as follows.
\begin{theorem}\label{theorem:derivatives} 
Assume that $f \in C^{1,2}\left([0,T]\times \mathbb{R}^n\right)$ and $\tbi{\E}{t_0}{x_0}\left[\left|	\mathcal{L}_{t, X_t} f(t, X_t)\right|\right]<\infty,$ let $t_0<t$ be a fixed time, and $x_0\in \mathbb{R}^n$ be a fixed space point, it holds that
\begin{equation*}
	\frac{d\,\tbi{\E}{t_0}{x_0}[f(t, X_t)]}{d t}=\tbi{\E}{t_0}{x_0}\left[\mathcal{A}_{t}^{ X_t} f(t, X_t)\right],\quad t\geq t_0.
\end{equation*}
Furthermore, one has the following identity
\begin{equation}\label{eq:PartialTimeX}
\left.\frac{d\,\tbi{\E}{t_0}{x_0}[f(t, X_t)]}{d t}\right|_{t=t_0}=\left.\frac{d\,\tbi{\E}{t_0}{x_0}[f(t, \tilde{X}_t)]}{d t}\right|_{t=t_0},
\end{equation}
where $\tilde{X}_t$ is an approximating diffusion process defined by
\begin{equation*}
\tilde{X}_t = x_0 + \int_0^t \tilde{a}\,ds + \int_0^t \tilde{b}\,dW_s
\end{equation*}
with the smooth functions $\tilde{a}_t=\tilde{a}(t,\tilde{X}_t;t_0,x_0)$ and $\tilde{b}_t=\tilde{b}(t,\tilde{X}_t;t_0,x_0)$ of $(t,\tilde{X}_t)$ with the parameter $(t_0,x_0)$ satisfying
\begin{equation*}
\tilde{a}(t_0,\tilde{X}_{t_0};t_0,x_0)=a(t_0,x_0)~\mbox{and}~\tilde{b}(t_0,\tilde{X}_{t_0};t_0,x_0)=b(t_0,x_0)
\end{equation*}
\end{theorem}
It has been noted in \cite{Zhao2014} that the different approximations of \eqref{eq:PartialTimeX} can be obtained by choosing different $\tilde{a}_t$'$s$ and $\tilde{b}_t$'$s.$ One can simply e.g., choose $ \tilde{a}(s,\tilde{X}_{s};t_0,x_0)=a(t_0,x_0)$ and $\tilde{b}(s,\tilde{X}_{s};t_0,x_0)=b(t_0,x_0)$ for all $s \in [t_0, t].$

For existence, regularity and representation of solutions of decoupled FBSDEs we refer to \cite{Ma2005, Peng1991, Zhang2001}. In the following of this section we will present some of those. We denote the forward stochastic differential equation (SDE) starting from $(s,x)$ with $X_t^{s,x}$ and consider the decoupled FBSDEs
\begin{equation}\label{eq:decoupledFBSDE}
\left\{
\begin{array}{l}
X_t^{s,x} = x + \int_s^t a(r, X^{s,x}_r)\,ds + \int_s^t b(r, X^{s,x}_r)\,dW_r,\\
Y^{s,x}_t = g(X^{s,x}_T) + \int_t^T f(r, X^{s,x}_r, Y^{s,x}_r, Z^{s,x}_r)\,dr - \int_t^T Z ^{s,x}_r\,dW_r,
\end{array}\right.
\end{equation} 
where $t\in [s, T],$ and the superscript $^{s,x}$ will be omitted when the context is clear.

Throughout the paper, we shall often make use of the following standing assumptions:
\begin{enumerate}
\item The functions $a, b \in C_b^1,$ and assume
\[
\sup_{0\leq t\leq T}\left\{|a(t,0)|+|b(t,0)|\right\} \leq L,
\]
where the common constant $L>0$ denotes all the Lipschitz constants.
\item $n=d$ and we assume that $b$ satisfies
\[b(t,x)b^\top(t,x) \geq \frac{1}{L} I_n,\quad \forall (t,x)\in [0,T]\times\mathbb{R}^n.\]
\item $a, b, f, g\in C_L,$ and assume that
\[
\sup_{0\leq t\leq T}|f(t,0,0,0)|+|g(0)| \leq L,
\]
where $L$ denotes all the Lipschitz constants.
\item $a, b, f \in C_L^{\frac{1}{2}}.$
\end{enumerate}
Under the above conditions, it is clear that \eqref{eq:decoupledFBSDE} is well-posed; the resulting integrands by taking conditional expectation on both side of the backward component is continuous with respect to time; the nonlinear Feynman-Kac formula \cite{Ma2005, Peng1991} can be given as follows. 
\begin{theorem}\label{theorem:FeynmanKac}
Let $u \in C^{1,2}\left([0,T]\times\mathbb{R}^n\right)$ be a classical solution to the following PDE
\[\mathcal{L}_{t,x}u(t,x)+f(t,x,u(t,x),\nabla u(t,x)b(t,x))=0,\quad u(T,x)=g(x),\]
then $ Y^{s,x}_t = u(t, X_t^{s,x}),$ $Z^{s,x}_t=\nabla_x u(t, X_t^{s,x})b(t,X_t^{s,x}),\, \forall t\in(s,T]$ is the unique solution to \eqref{eq:decoupledFBSDE}.
\end{theorem}
\section{Calculation of the weights in the FDM for approximating derivative}
In this section we calculate the weights in the FDM for approximating the function derivatives, e.g., $\frac{du(t)}{dt}.$
Let $u(t)\in C_b^{k+1}, k$ is a positive integer, and $t_i=i\Delta t,$ i.e., $t_0<t_1<\cdots<t_k.$
\subsection{Combination of two time points}\label{sub:twostep}
We consider the Taylor's expansions of $u(t_i)$ and $u(t_{i+1}), i=0,\cdots,k$
\begin{align*}
\begin{cases}
u(t_i)=\displaystyle\sum_{j=0}^{k}\frac{(\Delta t_i)^j}{j!}\frac{d^ju}{dt^j}(t_0)+\mathcal{O}(\Delta t_i)^{k+1}\\
u(t_{i+1})=\displaystyle\sum_{j=0}^{k}\frac{(\Delta t_{i+1})^j}{j!}\frac{d^ju}{dt^j}(t_0)+\mathcal{O}(\Delta t_{i+1})^{k+1},
\end{cases}
\end{align*}
from which we can deduce
\begin{align*}
\begin{cases}
\displaystyle\sum_{i=0}^{k}\alpha_{k,i} u(t_i)=\displaystyle\sum_{j=0}^{k}\frac{\displaystyle\sum_{i=0}^{k}\alpha_{k,i}(\Delta t_i)^j}{j!}\frac{d^ju}{dt^j}(t_0)+\mathcal{O}\left(\displaystyle\sum_{i=0}^{k}\alpha_{k,i}(\Delta t_i)^{k+1}\right)\\
\displaystyle\sum_{i=0}^{k}\alpha_{k,i}u(t_{i+1})=\displaystyle\sum_{j=0}^{k}\frac{\displaystyle\sum_{i=0}^{k}\alpha_{k,i}(\Delta t_{i+1})^j}{j!}\frac{d^ju}{dt^j}(t_0)+\mathcal{O}\left(\displaystyle\sum_{i=0}^{k}\alpha_{k,i}(\Delta t_{i+1})^{k+1}\right),
\end{cases}
\end{align*}
where $\alpha_{k,i}, i=0,1,\cdots,k$ are real numbers.
Clearly, we obtain
\begin{equation*}
\begin{split}
	\displaystyle\sum_{i=0}^{k}\alpha_{k,i} (u(t_i)+u(t_{i+1}))=\displaystyle\sum_{j=0}^{k}&\frac{\displaystyle\sum_{i=0}^{k}\alpha_{k,i}\left((\Delta t_i)^j+(\Delta t_{i+1})^j\right)}{j!}\frac{d^ju}{dt^j}(t_0)\\
	&\qquad \qquad +\mathcal{O}\left(\displaystyle\sum_{i=0}^{k}\alpha_{k,i}\left((\Delta t_i)^{k+1}+(\Delta t_{i+1})^{k+1}\right)\right)
\end{split}
\end{equation*}
and thus
\begin{equation}\label{eq:ODEApproximation}
	\frac{du}{dt}(t_0)=\displaystyle\sum_{i=0}^{k}\alpha_{k,i}\left(u(t_i)+u(t_{i+1})\right)+ \mathcal{O}\left(\displaystyle\sum_{i=0}^{k}\alpha_{k,i}\left((\Delta t_i)^{k+1}+(\Delta t_{i+1})^{k+1}\right)\right)
\end{equation}
by choosing
\begin{equation}\label{eq:systemCondition}
\frac{\displaystyle\sum_{i=0}^{k}\alpha_{k,i}\left((\Delta t_i)^j+(\Delta t_{i+1})^j\right)}{j!}=
\begin{cases}
1,\quad j=1,\\
0, \quad j \neq 1.
\end{cases}
\end{equation}
Due to $t_i=i\Delta t$ and $t_{i+1}=(i+1)\Delta t$ the conditions in \eqref{eq:systemCondition} are equivalent to the following system:
\[
\begin{bmatrix}
2 & 2 & 2 & \dots & 2 \\
1 & 3 & 5 & \dots & k+(k+1) \\
1 & 5 & 13 & \dots & k^2+(k+1)^2 \\
\vdots  & \vdots  & \vdots  & \vdots & \vdots  \\
1 & 1^k+2^k & 2^k+3^k & \dots & k^k+(k+1)^k 
\end{bmatrix}
\times
\begin{bmatrix}
\alpha_{k,0}\Delta t \\ \alpha_{k,1}\Delta t \\\alpha_{k,2}\Delta t \\\vdots \\ \alpha_{k,k}\Delta t 
\end{bmatrix}
=
\begin{bmatrix}
0 \\ 1 \\0\\ \vdots \\ 0
\end{bmatrix}
\]
which can be solved for $\alpha_{k,i}\Delta t, i=0,\cdots,k.$ We refer to the algorithm proposed in \cite{Fornberg1988} for those solutions. We report $\alpha_{k,i}\Delta t$ for $k=1, 2,\cdots, 7$ in Table \ref{tab:1}, since the related multi-step schemes proposed in this paper is unstable from $k=8,$ which will be explained below. 
\begin{table}[h]
	\centering
	\begin{tabular}{|c|c|c|c|c|c|c|c|c|}
		\hline
$\alpha_{k,i}\Delta t$&$i=0$ & $i=1$ &$i=2$&$i=3$ &$i=4$ &$i=5$ &$i=6$& $i=7$\\ \hline
$k=1$     & $-\frac{1}{2}$  & $\frac{1}{2}$ & & & & & & \\    \hline
$k=2$   & $-1$  & $\frac{3}{2}$ &$-\frac{1}{2}$ & & & & &\\    \hline
$k=3$     & $-\frac{17}{12}$  & $\frac{11}{4}$ &$-\frac{7}{4}$ &$\frac{5}{12}$ & & & &\\    \hline
$k=4$     & $-\frac{7}{4}$  & $\frac{49}{12}$ &$-\frac{15}{4}$ &$\frac{7}{4}$ & $-\frac{1}{3}$ & & & \\    \hline
$k=5$     & $-\frac{121}{60}$  & $\frac{65}{12}$ &$-\frac{77}{12}$ &$\frac{53}{12}$ &$-\frac{5}{3}$ & $\frac{4}{15}$& &\\    \hline
$k=6$     & $-\frac{67}{30}$  & $\frac{403}{60}$ &$-\frac{29}{3}$ &$\frac{35}{4}$ &$-\frac{59}{12}$ &$\frac{47}{30}$ &$-\frac{13}{60}$&\\    \hline
$k=7$     & $-\frac{2027}{840}$  & $\frac{319}{40}$ &$-\frac{1613}{120}$ &$\frac{361}{24}$ &$-\frac{269}{24}$ &$\frac{641}{120}$ &$-\frac{59}{40}$&$\frac{151}{840}$\\    \hline
	\end{tabular}
	\caption{The values of $\alpha_{k,i}\Delta t$ for combining two time points.}
	\label{tab:1}
\end{table}
The multi-step schemes (combining two time points) can be constructed by approximating the reference ODEs (see Sec. \ref{sec:TwoODEs}) using \eqref{eq:ODEApproximation}. Therefore, we consider the following ODE
\begin{equation}\label{eq:ODE}
	\frac{Y(t)}{dt}=f(t,Y(t)),\quad t\in[0,T)
\end{equation}
with the known terminal condition $Y(T)$ for studying stability, see also \cite{Zhao2014}. Applying \eqref{eq:ODEApproximation} to \eqref{eq:ODE} one obtain the multi-step scheme as 
\begin{equation}\label{eq:ODEMultistep}
	\alpha_{k,0}Y^n + \displaystyle\sum_{j=1}^{k}\left(\alpha_{k,j-1}+\alpha_{k,j}\right)Y^{n+j}+\alpha_{k,k}Y^{n+k+1}=f(t_n, Y^n)
\end{equation}
under the uniform time partition $0=t_0<t_1<\cdots<t_N=T.$ \eqref{eq:ODEMultistep} is stable if the roots $\{\lambda_{k,j}\}^k_{j=1}$ of the characteristic equation 
\begin{equation}\label{eq:CharateristicFun}
P(\lambda)=\alpha_{k,0} \lambda^{k+1} + \displaystyle\sum_{j=1}^{k}\left(\alpha_{k,j-1}+\alpha_{k,j}\right)\lambda^{k+1-j} + \alpha_{k,k}\lambda^0
\end{equation}
satisfies the following root conditions \cite{Butcher2008}
\begin{itemize}
\item $|\lambda_{k,j}|\leq 1,$
\item $P^{'}(\lambda_{k,j})\neq 0$ if $|\lambda_{k,j}|=1$ (simple roots).
\end{itemize}
With the $\alpha_{k,j}$ in Tabel \ref{tab:1}, $1$ is the simple root of the latter characteristic function for each $k,$ except which we list the maximum absolute values of the roots for $k=2,\cdots,8$ in Table \ref{tab:2}, from which we see that the multi-step scheme \eqref{eq:ODEMultistep} is unstable for $k\geq 8.$
\begin{table}[h]
	\centering
	\begin{tabular}{|c|c|c|c|c|c|c|c|}
		\hline
		$k$&$2$ & $3$ &$4$&$5$ &$6$ &$7$ &$8$\\ \hline
		$\max\left(\left|\lambda_{k,j}\right|\right)$     & $0.5000$  & $0.5424$ &$0.6344$ &$0.7438$ &$0.8636$ &$0.9915$ &  $1.1264$ \\   \hline
	\end{tabular}
	\caption{The maximum absolute root of \eqref{eq:CharateristicFun} except the simple roots}
	\label{tab:2}
\end{table}
However, compared to the multi-step scheme proposed in \cite{Zhao2014}(unstable $\geq 7$), stability for $k=7$ has been achieved, i.e., 1-order higher convergence rate is obtained. Combination of more time points can be done similarly, and provide other multi-step schemes, which have different instabilities. In our investigation we find that the multi-step scheme resulted by combining four time points are stable for $k\leq 9,$ which is the best. Thus, we show its detailed derivation in next subsection and will consider it in the numerical experiments.  
\subsection{Combination of four time points}
Similarly but slightly different to the multi-step scheme in Section \ref{sub:twostep}, we need to consider the Taylor's expansions of $u(t_i),$ $u(t_{i+1}),$ $u(t_{i+2})$ and $u(t_{i+3}), i=0,\cdots,k,$
\begin{align*}
\begin{cases}
u(t_i)=\displaystyle\sum_{j=0}^{k}\frac{(\Delta t_i)^j}{j!}\frac{d^ju}{dt^j}(t_0)+\mathcal{O}(\Delta t_i)^{k+1},\\
u(t_{i+1})=\displaystyle\sum_{j=0}^{k}\frac{(\Delta t_{i+1})^j}{j!}\frac{d^ju}{dt^j}(t_0)+\mathcal{O}(\Delta t_{i+1})^{k+1},\\
u(t_{i+2})=\displaystyle\sum_{j=0}^{k}\frac{(\Delta t_{i+2})^j}{j!}\frac{d^ju}{dt^j}(t_0)+\mathcal{O}(\Delta t_{i+2})^{k+1},\\
u(t_{i+3})=\displaystyle\sum_{j=0}^{k}\frac{(\Delta t_{i+3})^j}{j!}\frac{d^ju}{dt^j}(t_0)+\mathcal{O}(\Delta t_{i+3})^{k+1},
\end{cases}
\end{align*}
from which we can deduce
\begin{align*}
\begin{cases}
\displaystyle\sum_{i=0}^{k}\alpha_{k,i} u(t_i)=\displaystyle\sum_{j=0}^{k}\frac{\displaystyle\sum_{i=0}^{k}\alpha_{k,i}(\Delta t_i)^j}{j!}\frac{d^ju}{dt^j}(t_0)+\mathcal{O}\left(\sum_{i=0}^{k}\alpha_{k,i}(\Delta t_i)^{k+1}\right),\\
\displaystyle\sum_{i=0}^{k}\alpha_{k,i}u(t_{i+1})=\displaystyle\sum_{j=0}^{k}\frac{\displaystyle\sum_{i=0}^{k}\alpha_{k,i}(\Delta t_{i+1})^j}{j!}\frac{d^ju}{dt^j}(t_0)+\mathcal{O}\left(\displaystyle\sum_{i=0}^{k}\alpha_{k,i}(\Delta t_{i+1})^{k+1}\right),\\
\displaystyle\sum_{i=0}^{k}\alpha_{k,i}u(t_{i+2})=\displaystyle\sum_{j=0}^{k}\frac{\displaystyle\sum_{i=0}^{k}\alpha_{k,i}(\Delta t_{i+2})^j}{j!}\frac{d^ju}{dt^j}(t_0)+\mathcal{O}\left(\displaystyle\sum_{i=0}^{k}\alpha_{k,i}(\Delta t_{i+2})^{k+1}\right),\\
\displaystyle\sum_{i=0}^{k}\alpha_{k,i}u(t_{i+3})=\displaystyle\sum_{j=0}^{k}\frac{\displaystyle\sum_{i=0}^{k}\alpha_{k,i}(\Delta t_{i+3})^j}{j!}\frac{d^ju}{dt^j}(t_0)+\mathcal{O}\left(\displaystyle\sum_{i=0}^{k}\alpha_{k,i}(\Delta t_{i+3})^{k+1}\right),
\end{cases}
\end{align*}
where $\alpha_{k,i}, i=0,1,\cdots,k$ are real numbers as well.
Straightforwardly, we obtain
\begin{equation}\label{eq:errorTerm}
\begin{split}
\displaystyle\sum_{i=0}^{k}\alpha_{k,i} (u(t_i)&+u(t_{i+1})+u(t_{i+2})+u(t_{i+3}))=\\
&\displaystyle\sum_{j=0}^{k}\frac{\displaystyle\sum_{i=0}^{k}\alpha_{k,i}\left((\Delta t_i)^j+(\Delta t_{i+1})^j+(\Delta t_{i+2})^j+(\Delta t_{i+3})^j\right)}{j!}\frac{d^ju}{dt^j}(t_0)\\
&+\underbrace{\mathcal{O}\left(\displaystyle\sum_{i=0}^{k}\alpha_{k,i}\left((\Delta t_i)^{k+1}+(\Delta t_{i+1})^{k+1}+(\Delta t_{i+2})^{k+1}+(\Delta t_{i+3})^{k+1}\right)\right)}_{:=\epsilon}
\end{split}
\end{equation}
and thus
\begin{equation}\label{eq:multistepscheme}
\frac{du}{dt}(t_0)=\displaystyle\sum_{i=0}^{k}\alpha_{k,i}\left(u(t_i)+u(t_{i+1})+u(t_{i+2})+u(t_{i+3})\right)+\epsilon
\end{equation}
by choosing
\begin{equation*}
\frac{\displaystyle\sum_{i=0}^{k}\alpha_{k,i}\left((\Delta t_i)^j+(\Delta t_{i+1})^j+(\Delta t_{i+2})^j+(\Delta t_{i+3})^j\right)}{j!}=
\begin{cases}
1,\quad j=1,\\
0, \quad j \neq 1
\end{cases}
\end{equation*}
which are equivalent to the following system:
\[
\begin{bmatrix}
4 & 4 &  & \dots & 4 \\
6 & 10 & & \dots & k+(k+1)+(k+2)+(k+3) \\
14 & 30 &  & \dots & k^2+(k+1)^2+(k+2)^2+(k+3)^2\\
\vdots  & \vdots  &   & \vdots & \vdots  \\
1^k+2^k+3^k & 1^k+2^k+3^k+4^k &  & \dots & k^k+(k+1)^k+(k+2)^k+(k+3)^k
\end{bmatrix}
\times
\begin{bmatrix}
\alpha_{k,0}\Delta t \\ \alpha_{k,1}\Delta t \\\alpha_{k,2}\Delta t \\\vdots \\ \alpha_{k,k}\Delta t 
\end{bmatrix}
=
\begin{bmatrix}
0 \\ 1 \\0\\ \vdots \\ 0
\end{bmatrix}.
\]
In Table \ref{tab:3} we report solutions of the latter system for $k=1,\cdots, 9.$ 
\begin{table}[h]
	\centering
	\begin{tabular}{|c|c|c|c|c|c|c|c|c|c|c|}
		\hline
		$\alpha_{k,i}\Delta t$&$i=0$ & $i=1$ &$i=2$&$i=3$ &$i=4$ &$i=5$ &$i=6$& $i=7$&$i=8$&$i=9$\\ \hline
		$k=1$     & $-\frac{1}{4}$  & $\frac{1}{4}$ & & & & & & & &  \\    \hline
		$k=2$   & $-\frac{3}{4}$  & $\frac{5}{4}$ &$-\frac{1}{2}$ & & & & & & &\\    \hline
		$k=3$     & $-\frac{4}{3}$  & $3$ &$-\frac{9}{4}$ &$\frac{7}{12}$ & & & && &\\    \hline
		$k=4$     & $-\frac{11}{6}$  & $5$ &$-\frac{21}{4}$ &$\frac{31}{12}$ & $-\frac{1}{2}$ & & & & & \\    \hline
		$k=5$     & $-\frac{87}{40}$  & $\frac{161}{24}$ &$-\frac{26}{3}$ &$6$ &$-\frac{53}{24}$ & $\frac{41}{120}$& & & &\\    \hline
		$k=6$     & $-\frac{19}{8}$  & $\frac{949}{120}$ &$-\frac{35}{3}$ &$10$ &$-\frac{125}{24}$ &$\frac{37}{24}$ &$-\frac{1}{5}$& & &\\    \hline
		$k=7$     & $-\frac{419}{168}$  & $\frac{1049}{120}$ &$-\frac{85}{6}$ &$\frac{85}{6}$ &$-\frac{75}{8}$ &$\frac{97}{24}$ &$-\frac{31}{30}$&$\frac{5}{42}$& &\\    \hline
			$k=8$     & $-\frac{145}{56}$  & $\frac{2661}{280}$ &$-\frac{101}{6}$ &$\frac{39}{2}$ &$-\frac{385}{24}$ &$\frac{75}{8}$ &$-\frac{37}{10}$&$\frac{37}{42}$& -$\frac{2}{21}$&\\    \hline
				$k=9$     & $-\frac{6781}{2520}$  & $\frac{2917}{280}$ &$-\frac{4303}{210}$ &$\frac{841}{30}$ &$-\frac{3461}{120}$ &$\frac{887}{40}$ &$-\frac{367}{30}$&$\frac{953}{210}$&-$\frac{106}{105}$ &$\frac{32}{315}$\\    \hline
	\end{tabular}
	\caption{The values of $\alpha_{k,i}\Delta t$ for combining four time points.}
	\label{tab:3}
\end{table}
Applying \eqref{eq:multistepscheme} to \eqref{eq:ODE} one obtain the multi-step scheme as 
\begin{equation}\label{eq:MainODEMultistep}
\begin{split}
\alpha_{k,0}Y^n &+(\alpha_{k,0}+\alpha_{k,1})Y^{n+1}+(\alpha_{k,0}+\alpha_{k,1}+\alpha_{k,2})Y^{n+2}\\
&+\displaystyle\sum_{j=3}^{k}\left(\alpha_{k,j-3}+\alpha_{k,j-2}+\alpha_{k,j-1}+\alpha_{k,j}\right)Y^{n+j}+(\alpha_{k,k-2}+\alpha_{k,k-1}+\alpha_{k,k})Y^{n+k+1}
\\&+(\alpha_{k,k-1}+\alpha_{k,k})Y^{n+k+2}+\alpha_{k,k}Y^{n+k+3}=f(t_n, Y^n)
\end{split}
\end{equation}
whose characteristic equation reads
\begin{equation}\label{eq:MainCharateristicFun}
\begin{split}
	\alpha_{k,0}\lambda^{k+3} &+(\alpha_{k,0}+\alpha_{k,1})\lambda^{k+2}+(\alpha_{k,0}+\alpha_{k,1}+\alpha_{k,2})\lambda^{k+1}\\
	&+\displaystyle\sum_{j=3}^{k}\left(\alpha_{k,j-3}+\alpha_{k,j-2}+\alpha_{k,j-1}+\alpha_{k,j}\right)\lambda^{k+3-j}
	\\&+(\alpha_{k,k-2}+\alpha_{k,k-1}+\alpha_{k,k})\lambda^{2}+(\alpha_{k,k-1}+\alpha_{k,k})\lambda^{1}+\alpha_{k,k}\lambda^{0}=0.
\end{split}
\end{equation}
With the $\alpha_{k,j}$ in Table \ref{tab:3}, $1$ is the simple root of the latter characteristic function for each $k.$ The maximum absolute values of the roots for $k=2,\cdots,10$ expect the simple roots are listed in Table \ref{tab:4}, also the multi-step scheme \eqref{eq:MainODEMultistep} is stable for $k\leq 9.$
\begin{table}[h]
	\centering
	\begin{tabular}{|c|c|c|c|c|c|c|c|c|c|}
		\hline
		$k$&$2$ & $3$ &$4$&$5$ &$6$ &$7$ &$8$&$9$&$10$\\ \hline
		$\max\left(\left|\lambda_{k,j}\right|\right)$     & $0.6667$ & $0.6614$ &$0.6875$ &$0.7104$ &$0.7224$ &$0.7376$ &$0.8134$  &$0.9931$& $1.2286$\\    \hline
	\end{tabular}
	\caption{The maximum absolute root of \eqref{eq:MainCharateristicFun} except the simple roots}
	\label{tab:4}
\end{table}
We remark that the stability cannot be guaranteed for $k>9$ by combining more time points, e.g., the multi-step scheme constructed by combining five time points is stable for $k\leq 8.$
\section{The semi-discrete multi-step scheme for decoupled FBSDEs}
Following the idea in \cite{Zhao2014} we derive the semi-discrete scheme for \eqref{eq:fbsde} in the decoupled case. We consider the time interval $[0,T]$ with the following partition
\begin{equation*}
	0=t_0<t_1<t_2<\cdots t_{N_T}=T.
\end{equation*}
We denote $t_{n+k}-t_n$ by $\Delta t_{n,k}$ and $W_{t_{n+k}} -W_{t_n}$ by $\Delta W_{n,k},$ i.e., $\Delta t_{t_n,t}=t-t_n$ and $\Delta W_{t_n,t}=W_t-W_{t_n}$ for $t\geq t_n.$
\subsection{Two reference ODEs}\label{sec:TwoODEs}
Let $(X_t,Y_t,Z_t)$ be the solution of the decoupled FBSDEs \eqref{eq:fbsde}. By taking conditional expectation $\E_{t_n}^x[\cdot]$ on both sides of the backward component in \eqref{eq:fbsde} one obtains the integral equation
\begin{equation*}
\E_{t_n}^x\left[Y_t\right]=\E_{t_n}^x\left[\xi\right] + \int_t^T\E_{t_n}^x\left[f(s,X_s,Y_s,Z_s)\right]\,ds,\quad \forall t \in [t_n,T].
\end{equation*}
As explained in Sec. \ref{sec:Preliminaries}, the integrand in the latter integral equation is continuous with respect to the time. By taking the derivative with respect to $t$ on both sides one thus obtain the first reference ODE:
\begin{equation}\label{eq:FirstODE}
\frac{d\,\E_{t_n}^x\left[Y_t\right]}{dt}=-\E_{t_n}^x\left[f(t,X_t,Y_t,Z_t)\right],\quad \forall t \in [t_n,T].
\end{equation}
Furthermore, we have
\begin{equation*}
	Y_{t_n}=Y_t+\int_{t_n}^{t}f(s,X_s,Y_s,Z_s)\,ds-\int_{t_n}^{t}Z_s\,dW_s,\quad t\in[t_n,T].
\end{equation*}
By multiplying both sides of the latter equation by $\Delta W_{t_n,t}^{\top}$ and again taking the conditional expectation $\E_{t_n}^x[\cdot]$ on its both sides we obtain
\begin{equation*}
	0=\E_{t_n}^x\left[Y_t\Delta W_{t_n,t}^{\top}\right]+\int_{t_n}^{t}\E_{t_n}^x\left[f(s,X_s,Y_s,Z_s)\Delta W_{t_n,s}^{\top}\right]\,ds-\int_{t_n}^{t}\E_{t_n}^x\left[Z_s\right]\,ds,\quad t\in[t_n,T].
\end{equation*}
Similarly, we obtain the second reference ODE:
\begin{equation}\label{eq:SecondODE}
\frac{d\,\E_{t_n}^x\left[Y_t\Delta W_{t_n,t}^{\top}\right]}{dt}=-\E_{t_n}^x\left[f(t,X_t,Y_t,Z_t)\Delta W_{t_n,t}^{\top}\right]+ \E_{t_n}^x\left[Z_t\right],\quad t \in [t_n,T].
\end{equation}
by taking the derivative with respect to $t \in [t_n, T].$
\subsection{The semi-discrete scheme}
Let $\bar{a}(t,x)$ and $\bar{b}(t,x)$ be smooth functions for $t\in[t_n,T]$ and $x\in \mathbb{R}^n$ satisfying $\bar{b}(t,x)=b(t,x)$ and $\bar{b}(t,x)=b(t,x),$ and thus define the diffusion process
\begin{equation}\label{eq:Xapproximation}
	\bar{X}_{t}^{t_n,x} = x + \int_{t_n}^{t}\bar{b}(s,\bar{X}_{s}^{t_n,x})\,ds+\int_{t_n}^{t}\bar{b}(s,\bar{X}_{s}^{t_n,x})\,dW_s.
\end{equation}
Let $(X_{t}^{t_n,x},Y_{t}^{t_n,x},Z_{t}^{t_n,x})$ be the solution of the decoupled FBSDEs, i.e., $Y_{t}^{t_n,x}$ and $Z_{t}^{t_n,x}$ can be represented by $u(t,X_{t}^{t_n,x})$ and $\nabla_x u(t,X_{t}^{t_n,x}) b(s,X_{s}^{t_n,x}),$ respectively, see Theorem \ref{theorem:FeynmanKac}.

Therefore, we set $\bar{Y}_t^{t_n,x}=u(t,\bar{X}_{t}^{t_n,x})$ and $\bar{Z}_t^{t_n,x}=\nabla_x u(t,\bar{X}_{t}^{t_n,x}) b(s,\bar{X}_{s}^{t_n,x})$ to have
\begin{equation*}
		\left.\frac{d\,\tbi{\E}{t_n}{x}[\tbi{Y}{t}{t_n,x}]}{d t}\right|_{t=t_n}=		\left.\frac{d\,\tbi{\E}{t_n}{x}[\tbi{\bar{Y}}{t}{t_n,x}]}{d t}\right|_{t=t_n}
\end{equation*}
and
\begin{equation*}
	\left.\frac{d\,\tbi{\E}{t_n}{x}[\tbi{Y}{t}{t_n,x}\Delta W_{t_n,t}^{\top}]}{d t}\right|_{t=t_n}=		\left.\frac{d\,\tbi{\E}{t_n}{x}[\tbi{\bar{Y}}{t}{t_n,x}\Delta W_{t_n,t}^{\top}]}{d t}\right|_{t=t_n}
\end{equation*}
by Theorem \ref{theorem:derivatives}. Then, we apply \eqref{eq:multistepscheme} to terms on the right hand side of both the latter equations to obtain
\begin{equation}\label{eq:ForODE1}
	\begin{split}
	\left.\frac{d\,\tbi{\E}{t_n}{x}[\tbi{Y}{t}{t_n,x}]}{d t}\right|_{t=t_n} &= \displaystyle\sum_{i=0}^{k}\alpha_{k,i}\tbi{\E}{t_n}{x}\left[\tbi{\bar{Y}}{t_{n+i}}{t_n,x}+\tbi{\bar{Y}}{t_{n+i+1}}{t_n,x}+\tbi{\bar{Y}}{t_{n+i+2}}{t_n,x}+\tbi{\bar{Y}}{t_{n+i+3}}{t_n,x}\right]+\tbi{\bar{R}}{y,n}{k}\\
	=\alpha_{k,0}\tbi{\E}{t_n}{x}\left[\tbi{\bar{Y}}{t_{n}}{t_n,x}\right] &+(\alpha_{k,0}+\alpha_{k,1})\tbi{\E}{t_n}{x}\left[\tbi{\bar{Y}}{t_{n+1}}{t_n,x}\right]+(\alpha_{k,0}+\alpha_{k,1}+\alpha_{k,2})\tbi{\E}{t_n}{x}\left[\tbi{\bar{Y}}{t_{n+2}}{t_n,x}\right]\\
	&+\displaystyle\sum_{j=3}^{k}\left(\alpha_{k,j-3}+\alpha_{k,j-2}+\alpha_{k,j-1}+\alpha_{k,j}\right)\tbi{\E}{t_n}{x}\left[\tbi{\bar{Y}}{t_{n+j}}{t_n,x}\right]
	\\&+(\alpha_{k,k-2}+\alpha_{k,k-1}+\alpha_{k,k})\tbi{\E}{t_n}{x}\left[\tbi{\bar{Y}}{t_{n+k+1}}{t_n,x}\right]\\
	&+(\alpha_{k,k-1}+\alpha_{k,k})\tbi{\E}{t_n}{x}\left[\tbi{\bar{Y}}{t_{n+k+2}}{t_n,x}\right]+\alpha_{k,k}\tbi{\E}{t_n}{x}\left[\tbi{\bar{Y}}{t_{n+k+3}}{t_n,x}\right]+\tbi{\bar{R}}{y,n}{k}
	\end{split}
\end{equation}
and
\begin{equation}\label{eq:ForODE2}
\begin{split}
&\left.\frac{d\,\tbi{\E}{t_n}{x}[\tbi{Y}{t}{t_n,x}\Delta W_{t_n,t}^{\top}]}{d t}\right|_{t=t_n}=\displaystyle\sum_{i=1}^{k}\alpha_{k,i}\tbi{\E}{t_n}{x}\left[\tbi{\bar{Y}}{t_{n+i}}{t_n,x}\Delta W_{n,i}^{\top}+\tbi{\bar{Y}}{t_{n+i+1}}{t_n,x}\Delta W_{n,i+1}^{\top}\right.\\
&\qquad\left.+\tbi{\bar{Y}}{t_{n+i+2}}{t_n,x}\Delta W_{n,i+2}^{\top}+\tbi{\bar{Y}}{t_{n+i+3}}{t_n,x}\Delta W_{n,i+3}^{\top}\right]+\tbi{\bar{R}}{z,n}{k}\\
&=(\alpha_{k,0}+\alpha_{k,1})\tbi{\E}{t_n}{x}\left[\tbi{\bar{Y}}{t_{n+1}}{t_n,x}\Delta W_{n,1}^{\top}\right]+(\alpha_{k,0}+\alpha_{k,1}+\alpha_{k,2})\tbi{\E}{t_n}{x}\left[\tbi{\bar{Y}}{t_{n+2}}{t_n,x}\Delta W_{n,2}^{\top}\right]\\
&+\displaystyle\sum_{j=3}^{k}\left(\alpha_{k,j-3}+\alpha_{k,j-2}+\alpha_{k,j-1}+\alpha_{k,j}\right)\tbi{\E}{t_n}{x}\left[\tbi{\bar{Y}}{t_{n+j}}{t_n,x}\Delta W_{n,j}^{\top}\right]
\\&+(\alpha_{k,k-2}+\alpha_{k,k-1}+\alpha_{k,k})\tbi{\E}{t_n}{x}\left[\tbi{\bar{Y}}{t_{n+k+1}}{t_n,x}\Delta W_{n,k+1}^{\top}\right]\\
&+(\alpha_{k,k-1}+\alpha_{k,k})\tbi{\E}{t_n}{x}\left[\tbi{\bar{Y}}{t_{n+k+2}}{t_n,x}\Delta W_{n,k+2}^{\top}\right]+\alpha_{k,k}\tbi{\E}{t_n}{x}\left[\tbi{\bar{Y}}{t_{n+k+3}}{t_n,x}\Delta W_{n,k+3}^{\top}\right]+\tbi{\bar{R}}{z,n}{k}
\end{split}
\end{equation}
where $\alpha_{k,i}$ are given in Table \ref{tab:3},
$\tbi{\bar{R}}{y,n}{k}$ and $\tbi{\bar{R}}{z,n}{k}$ are truncation errors.
We insert respectively \eqref{eq:ForODE1} and \eqref{eq:ForODE2} into \eqref{eq:FirstODE} and \eqref{eq:SecondODE}, and obtain
\begin{equation}\label{eq:BeforeScheme1}
\begin{split}
&\alpha_{k,0}\tbi{\E}{t_n}{x}\left[\tbi{\bar{Y}}{t_{n}}{t_n,x}\right] +(\alpha_{k,0}+\alpha_{k,1})\tbi{\E}{t_n}{x}\left[\tbi{\bar{Y}}{t_{n+1}}{t_n,x}\right]+(\alpha_{k,0}+\alpha_{k,1}+\alpha_{k,2})\tbi{\E}{t_n}{x}\left[\tbi{\bar{Y}}{t_{n+2}}{t_n,x}\right]\\
&+\displaystyle\sum_{j=3}^{k}\left(\alpha_{k,j-3}+\alpha_{k,j-2}+\alpha_{k,j-1}+\alpha_{k,j}\right)\tbi{\E}{t_n}{x}\left[\tbi{\bar{Y}}{t_{n+j}}{t_n,x}\right]
\\&+(\alpha_{k,k-2}+\alpha_{k,k-1}+\alpha_{k,k})\tbi{\E}{t_n}{x}\left[\tbi{\bar{Y}}{t_{n+k+1}}{t_n,x}\right]+(\alpha_{k,k-1}+\alpha_{k,k})\tbi{\E}{t_n}{x}\left[\tbi{\bar{Y}}{t_{n+k+2}}{t_n,x}\right]\\
&+\alpha_{k,k}\tbi{\E}{t_n}{x}\left[\tbi{\bar{Y}}{t_{n+k+3}}{t_n,x}\right]=-f(t_n,x,Y_{t_n},Z_{t_n})+\tbi{R}{y,n}{k}
\end{split}
\end{equation}
and
\begin{equation}\label{eq:BeforeScheme2}
\begin{split}
&(\alpha_{k,0}+\alpha_{k,1})\tbi{\E}{t_n}{x}\left[\tbi{\bar{Y}}{t_{n+1}}{t_n,x}\Delta W_{n,1}^{\top}\right]+(\alpha_{k,0}+\alpha_{k,1}+\alpha_{k,2})\tbi{\E}{t_n}{x}\left[\tbi{\bar{Y}}{t_{n+2}}{t_n,x}\Delta W_{n,2}^{\top}\right]\\
&+\displaystyle\sum_{j=3}^{k}\left(\alpha_{k,j-3}+\alpha_{k,j-2}+\alpha_{k,j-1}+\alpha_{k,j}\right)\tbi{\E}{t_n}{x}\left[\tbi{\bar{Y}}{t_{n+j}}{t_n,x}\Delta W_{n,j}^{\top}\right]
\\&+(\alpha_{k,k-2}+\alpha_{k,k-1}+\alpha_{k,k})\tbi{\E}{t_n}{x}\left[\tbi{\bar{Y}}{t_{n+k+1}}{t_n,x}\Delta W_{n,k+1}^{\top}\right]\\
&+(\alpha_{k,k-1}+\alpha_{k,k})\tbi{\E}{t_n}{x}\left[\tbi{\bar{Y}}{t_{n+k+2}}{t_n,x}\Delta W_{n,k+2}^{\top}\right]+\alpha_{k,k}\tbi{\E}{t_n}{x}\left[\tbi{\bar{Y}}{t_{n+k+3}}{t_n,x}\Delta W_{n,k+3}^{\top}\right]=Z_{t_n}+\tbi{R}{z,n}{k}
\end{split}
\end{equation}
with $\tbi{R}{y,n}{k}=-\tbi{\bar{R}}{y,n}{k}$ and $\tbi{R}{z,n}{k}=-\tbi{\bar{R}}{z,n}{k}.$

We denote the numerical approximations of $Y_t$ and $Z_t$ at $t_n$ by $Y^n$ and $Z^n,$ respectively. Furthermore, for $\bar{a}$ and $\bar{b}$ in \eqref{eq:Xapproximation} we choose $\bar{a}(t,\bar{X}_{t}^{t_n,x})=a(t_n,x)$ and $\bar{b}(t,\bar{X}_{t}^{t_n,x})=b(t_n,x)$ for $t\in [t_n,T].$ Finally, from \eqref{eq:BeforeScheme1} and \eqref{eq:BeforeScheme2}, the semi-discrete scheme can be obtained as
\begin{scheme}
Assume that $Y^{N_T-i}$ and $Z^{N_T-i}$ are known for $i=0,1,\cdots,k+2.$ For $n=N_T-k-3,\cdots, 0,$ $X^{n,j},$
$Y^{n}=Y^{n}(X^n)$ and $Z^{n}=Z^{n}(X^n)$ can be solved by
\begin{equation}\label{eq:SchemeX}
	\hspace{-2.2cm}X^{n,j}=X^n + a(t_n,X^n)\Delta t_{n,j}+ b(t_n,X^n)\Delta W_{n,j},\quad j=1,\cdots,k+3,
\end{equation}
\begin{equation}\label{eq:SchemeY}
\begin{split}
&Z^n=(\alpha_{k,0}+\alpha_{k,1})\tbi{\E}{t_n}{X^n}\left[\tbi{\bar{Y}}{}{n+1}\Delta W_{n,1}^{\top}\right]+(\alpha_{k,0}+\alpha_{k,1}+\alpha_{k,2})\tbi{\E}{t_n}{X^n}\left[\tbi{\bar{Y}}{}{n+2}\Delta W_{n,2}^{\top}\right]\\
&+\displaystyle\sum_{j=3}^{k}\left(\alpha_{k,j-3}+\alpha_{k,j-2}+\alpha_{k,j-1}+\alpha_{k,j}\right)\tbi{\E}{t_n}{X^n}\left[\tbi{\bar{Y}}{}{n+j}\Delta W_{n,j}^{\top}\right]
\\&+(\alpha_{k,k-2}+\alpha_{k,k-1}+\alpha_{k,k})\tbi{\E}{t_n}{X^n}\left[\tbi{\bar{Y}}{}{n+k+1}\Delta W_{n,k+1}^{\top}\right]\\
&+(\alpha_{k,k-1}+\alpha_{k,k})\tbi{\E}{t_n}{X^n}\left[\tbi{\bar{Y}}{}{n+k+2}\Delta W_{n,k+2}^{\top}\right]+\alpha_{k,k}\tbi{\E}{t_n}{X^n}\left[\tbi{\bar{Y}}{}{n+k+3}\Delta W_{n,k+3}^{\top}\right],
\end{split}
\end{equation}
\begin{equation}\label{eq:SchemeZ}
\begin{split}
\alpha_{k,0}Y^n&=-(\alpha_{k,0}+\alpha_{k,1})\tbi{\E}{t_n}{X^n}\left[\bar{Y}^{n+1}\right]-(\alpha_{k,0}+\alpha_{k,1}+\alpha_{k,2})\tbi{\E}{t_n}{X^n}\left[\bar{Y}^{n+2}\right]\\
&-\displaystyle\sum_{j=3}^{k}(\alpha_{k,j-3}+\alpha_{k,j-2}+\alpha_{k,j-1}+\alpha_{k,j})\tbi{\E}{t_n}{X^n}\left[\bar{Y}^{n+j}\right]\\
&-(\alpha_{k,k-2}+\alpha_{k,k-1}+\alpha_{k,k})\tbi{\E}{t_n}{X^n}\left[\bar{Y}^{n+k+1}\right]-(\alpha_{k,k-1}+\alpha_{k,k})\tbi{\E}{t_n}{X^n}\left[\bar{Y}^{n+k+2}\right]\\
&-\alpha_{k,k}\tbi{\E}{t_n}{X^n}\left[\bar{Y}^{n+k+3}\right]-f(t_n,X^n,Y^n,Z^n).
\end{split}	
\end{equation}
\end{scheme}
\begin{remark}
	\begin{enumerate}
	\item $\bar{Y}^{n+j}$ is the value of $Y^{n+j}$ at the space point $X^{n,j}$ for $j=1,\cdots,k+3.$
	\item The latter implicit equation can be solved by using iterative methods, e.g., Newton's method or Picard scheme. 
	\item By Theorem \ref{theorem:derivatives} and \eqref{eq:errorTerm} it holds \cite{Butcher2008}
	\begin{equation*}
		\tbi{\bar{R}}{y,n}{k}=\mathcal{O}(\Delta t)^k~\mbox{and}~\tbi{\bar{R}}{z,n}{k}=\mathcal{O}(\Delta t)^k
	\end{equation*}
	provided that $\mathcal{L}^{k+4}_{t,x}u(t,x)$ is bounded, where $\tbi{\bar{R}}{y,n}{k}$ and $\tbi{\bar{R}}{z,n}{k}$
	are defined in \eqref{eq:ForODE1} and \eqref{eq:ForODE2}, respectively.
	\item Similar to the scheme proposed in \cite{Zhao2014}, one can obtain high-order accurate numerical solutions for \eqref{eq:SchemeY} and \eqref{eq:SchemeZ}, although the Euler scheme is used for \eqref{eq:SchemeX}. This is the main advantages because the usage of the Euler scheme reduces
	dramatically the total computational complexity, and one is only interested in the solution of  \eqref{eq:SchemeY} and \eqref{eq:SchemeZ} in many applications.
	\end{enumerate}
\end{remark}
\section{The fully discrete multi-step scheme for decoupled FBSDEs}
To solve $(X^n, Y^n,Z^n)$ numerically, next we consider the space discretization. We define firstly the partition of the real space as
$\mathcal{R}^{n}_h=\{x_i|x_i\in \mathbb{R}^n\}$ with
\begin{equation*}
h^n=\max_{x \in \mathbb{R}^n} dist(x,\mathcal{R}_h^n),
\end{equation*}
where $dist(x,\mathcal{R}_h^n)$ is the distance from $x$ to $\mathcal{R}_h^x.$ Furthermore, for each $x$ we define the neighbor grid set (local subset)
$\mathcal{R}_{h,x}^n$ satisfying
\begin{enumerate}
\item $dist(x,\mathcal{R}_h^n) < dist(x,\mathcal{R}_h^n)/\ \mathcal{R}_{h,x}^n,$
\item the number of elements in $\mathcal{R}_{h,x}^n$ is finite and uniformly bounded.
\end{enumerate}
Based on the space discretization, we can solve $Y^n(x)$ and $Z^n(x)$ for each grid point $x\in \mathcal{R}^n_h,$
$n=N_t-k-3,\cdots,0,$ by
\begin{equation}\label{eq:FSchemeY}
\begin{split}
&Z^n=(\alpha_{k,0}+\alpha_{k,1})\tbi{\E}{t_n}{x}\left[\tbi{\bar{Y}}{}{n+1}\Delta W_{n,1}^{\top}\right]+(\alpha_{k,0}+\alpha_{k,1}+\alpha_{k,2})\tbi{\E}{t_n}{x}\left[\tbi{\bar{Y}}{}{n+2}\Delta W_{n,2}^{\top}\right]\\
&+\displaystyle\sum_{j=3}^{k}\left(\alpha_{k,j-3}+\alpha_{k,j-2}+\alpha_{k,j-1}+\alpha_{k,j}\right)\tbi{\E}{t_n}{x}\left[\tbi{\bar{Y}}{}{n+j}\Delta W_{n,j}^{\top}\right]
\\&+(\alpha_{k,k-2}+\alpha_{k,k-1}+\alpha_{k,k})\tbi{\E}{t_n}{x}\left[\tbi{\bar{Y}}{}{n+k+1}\Delta W_{n,k+1}^{\top}\right]\\
&+(\alpha_{k,k-1}+\alpha_{k,k})\tbi{\E}{t_n}{x}\left[\tbi{\bar{Y}}{}{n+k+2}\Delta W_{n,k+2}^{\top}\right]+\alpha_{k,k}\tbi{\E}{t_n}{x}\left[\tbi{\bar{Y}}{}{n+k+3}\Delta W_{n,k+3}^{\top}\right],
\end{split}
\end{equation}
\begin{equation}\label{eq:FSchemeZ}
\begin{split}
\alpha_{k,0}Y^n&=-(\alpha_{k,0}+\alpha_{k,1})\tbi{\E}{t_n}{x}\left[\bar{Y}^{n+1}\right]-(\alpha_{k,0}+\alpha_{k,1}+\alpha_{k,2})\tbi{\E}{t_n}{x}\left[\bar{Y}^{n+2}\right]\\
&-\displaystyle\sum_{j=3}^{k}(\alpha_{k,j-3}+\alpha_{k,j-2}+\alpha_{k,j-1}+\alpha_{k,j})\tbi{\E}{t_n}{x}\left[\bar{Y}^{n+j}\right]\\
&-(\alpha_{k,k-2}+\alpha_{k,k-1}+\alpha_{k,k})\tbi{\E}{t_n}{x}\left[\bar{Y}^{n+k+1}\right]-(\alpha_{k,k-1}+\alpha_{k,k})\tbi{\E}{t_n}{x}\left[\bar{Y}^{n+k+2}\right]\\
&-\alpha_{k,k}\tbi{\E}{t_n}{x}\left[\bar{Y}^{n+k+3}\right]-f(t_n,x,Y^n,Z^n).
\end{split}	
\end{equation}
Note that $\bar{Y}^{n+j}$ is the value of $Y^{n+j}$ at the space point $X^{n,j}$ generated by
\begin{equation*}
	X^{n,j}=x + a(t_n,x)\Delta t_{n,j}+ b(t_n,x)\Delta W_{n,j},\quad j=1,\cdots,k+3.
\end{equation*}
However, $X^{n,j}$ does not belong to $\mathcal{R}_h^{n+j}.$ This is to say that the value of $Y^{n+j}$ 
at $X^{n,j}$ needs to be approximated based on the values of $Y^{n+j}$ on $\mathcal{R}_h^{n+j},$ this can be done using a local interpolation.
By $LI^n_{h,X}F$ we denote the interpolated value of the function $F$ at space point $X\in \mathbb{R}^n$ by using the values of $F$
only in the neighbor grid set, namely $\mathcal{R}^n_{h,X}.$ Including the interpolations, \eqref{eq:FSchemeY} and \eqref{eq:FSchemeZ} become
\begin{equation}\label{eq:FSchemeYI}
\begin{split}
Z^n=(\alpha_{k,0}&+\alpha_{k,1})\tbi{\E}{t_n}{x}\left[LI^{n+1}_{h,X^{n,j}}\tbi{Y}{}{n+1}\Delta W_{n,1}^{\top}\right]\\
&+(\alpha_{k,0}+\alpha_{k,1}+\alpha_{k,2})\tbi{\E}{t_n}{x}\left[LI^{n+2}_{h,X^{n,j}}\tbi{Y}{}{n+2}\Delta W_{n,2}^{\top}\right]\\
&+\displaystyle\sum_{j=3}^{k}\left(\alpha_{k,j-3}+\alpha_{k,j-2}+\alpha_{k,j-1}+\alpha_{k,j}\right)\tbi{\E}{t_n}{x}\left[LI^{n+j}_{h,X^{n,j}}\tbi{Y}{}{n+j}\Delta W_{n,j}^{\top}\right]
\\&+(\alpha_{k,k-2}+\alpha_{k,k-1}+\alpha_{k,k})\tbi{\E}{t_n}{x}\left[\tbi{LI^{n+k+1}_{h,X^{n,j}}Y}{}{n+k+1}\Delta W_{n,k+1}^{\top}\right]\\
&+(\alpha_{k,k-1}+\alpha_{k,k})\tbi{\E}{t_n}{x}\left[LI^{n+k+2}_{h,X^{n,j}}\tbi{Y}{}{n+k+2}\Delta W_{n,k+2}^{\top}\right]\\
&+\alpha_{k,k}\tbi{\E}{t_n}{x}\left[\tbi{LI^{n+k+3}_{h,X^{n,j}}Y}{}{n+k+3}\Delta W_{n,k+3}^{\top}\right]+R_{z,n}^{k,LI_h},
\end{split}
\end{equation}
\begin{equation}\label{eq:FSchemeZI}
\begin{split}
\alpha_{k,0}Y^n=&-(\alpha_{k,0}+\alpha_{k,1})\tbi{\E}{t_n}{x}\left[LI^{n+1}_{h,X^{n,j}}Y^{n+1}\right]-(\alpha_{k,0}+\alpha_{k,1}+\alpha_{k,2})\tbi{\E}{t_n}{x}\left[LI^{n+2}_{h,X^{n,j}}Y^{n+2}\right]\\
&-\displaystyle\sum_{j=3}^{k}(\alpha_{k,j-3}+\alpha_{k,j-2}+\alpha_{k,j-1}+\alpha_{k,j})\tbi{\E}{t_n}{x}\left[LI^{n+j}_{h,X^{n,j}}Y^{n+j}\right]\\
&-(\alpha_{k,k-2}+\alpha_{k,k-1}+\alpha_{k,k})\tbi{\E}{t_n}{x}\left[LI^{n+k+1}_{h,X^{n,j}}Y^{n+k+1}\right]\\
&-(\alpha_{k,k-1}+\alpha_{k,k})\tbi{\E}{t_n}{x}\left[LI^{n+k+2}_{h,X^{n,j}}Y^{n+k+2}\right]\\
&-\alpha_{k,k}\tbi{\E}{t_n}{x}\left[LI^{n+k+3}_{h,X^{n,j}}Y^{n+k+3}\right]-f(t_n,x,Y^n,Z^n)+R_{y,n}^{k,LI_h}.
\end{split}	
\end{equation}

Furthermore, to approximate the conditional expectations in \eqref{eq:FSchemeYI} and \eqref{eq:FSchemeZI} we employ the Gauss-Hermite quadrature rule which
is an extension of the Gaussian quadrature method for approximating the value of integrals of the form $\int_{-\infty}^{\infty}\exp(-{\bf x}^2)g({\bf x})\,d{\bf x}$ by
\begin{equation}\label{eq:GaussApproximation}
\int_{-\infty}^{\infty}\dots\int_{-\infty}^{\infty}g({\bf x})\exp(-{\bf x}^{\top}{\bf x})\,d{\bf x} \approx \displaystyle\sum_{{\bf j}=1}^{L}w_{{\bf j}}g({\bf a_j}),
\end{equation}
where ${\bf x}=(x_1,\cdots,x_n)^{\top},$ ${\bf x}^{\top} {\bf x}= \displaystyle\sum_{j=1}^{n} x^2_j,$ $L$ is the number of used sample points,
${\bf j}=\left(j_1,j_2,\cdots,j_n\right),$ $\displaystyle\sum_{{\bf j=1}}^{L}=\displaystyle\sum_{j_1=1,\cdots,j_n=1}^{L,\cdots,L},$ ${\bf a_j}=(a_{j_1},\cdots,a_{j_n})$
and $\omega_{\bf j}=\displaystyle\prod_{i=1}^{n}\omega_{j_i},$ $\{a_{j_i}\}_{j_i=1}^L$ are the roots of the Hermite polynomial $H_L(x)$ of degree $L$
and $\{\omega_{j_i}\}_{j_i=1}^L$ are corresponding weights \cite{Abramowitz1972}. For a standard $n$-dimensional standard normal distributed random variable $X$
we know that
\begin{align*}
	\E\left[g(X)\right]&=\frac{1}{(2\pi)^{\frac{d}{2}}}\int_{-\infty}^{\infty}g({\bf x})\exp\left(-\frac{{\bf x}^{\top}{\bf x}}{2}\right)\,dx\\
	&=\frac{1}{(\pi)^{\frac{d}{2}}}\int_{-\infty}^{\infty}g(\sqrt{2}{\bf x})\exp\left(-{\bf x}^{\top}{\bf x}\right)\,dx\\
	&\stackrel{\ref{eq:GaussApproximation}}{=} \frac{1}{(\pi)^{\frac{d}{2}}}\displaystyle\sum_{{\bf j=1}}^{L}\omega_{{\bf j}}g({\bf a_j})+R^{GH}_L,
\end{align*}
where $R^{GH}_L$ is the truncation error of the Gauss-Hermite quadrature rule for $g.$

Now we consider the conditional expectations of the form $\tbi{\E}{t_n}{x}\left[LI^{n+j}_{h,X^{n,j}}\tbi{Y}{}{n+j}\Delta W_{n,j}^{\top}\right]$ and $\tbi{\E}{t_n}{x}\left[LI^{n+j}_{h,X^{n,j}}Y^{n+j}\right]$ in 
\eqref{eq:FSchemeYI} and \eqref{eq:FSchemeZI}. We know that $LI^{n+j}_{h,X^{n,j}}Y^{n+j}$ is the interpolated value of $\bar{Y}^{n+j},$ which is a function
of $X^{n,j}$ and can be represented by (Theorem \ref{theorem:FeynmanKac})
\begin{equation*}
\bar{Y}^{n+j}=Y^{n+j}\left(X^{n+j}\right)=Y^{n+j}\left(X^{n+j}+a(t_n,X^n)\Delta t_{n,j}+b(t_n,X^n)\Delta W_{n,j}\right),
\end{equation*}
with $\Delta W_{n,j} \sim \sqrt{\Delta t_{n,j}}N(0,I_n).$ Straightforwardly, we can approximate those conditional expectations as
\begin{equation*}
	\E^{{\bf x},h}_{t_n}\left[\bar{Y}^{n+j}\right]=\frac{1}{\pi^{\frac{d}{2}}}\sum_{{\bf j=1}}^{L}\omega_{{\bf j}}Y^{n+j}
	\left({\bf x}+a(t_n,{\bf x})\Delta t_{n,j}+b(t_n,{\bf x})\Delta t_{n,j}\sqrt{2\Delta t_{n,j}}{\bf a_j}\right)+R^{GH}_L(Y)
\end{equation*}
and
\begin{equation*}
\E^{{\bf x},h}_{t_n}\left[\bar{Y}^{n+j}\Delta W_{t_n, j}^{\top}\right]=\frac{1}{\pi^{\frac{d}{2}}}\sum_{{\bf j=1}}^{L}\omega_{{\bf j}}Y^{n+j}
\left({\bf x}+a(t_n,{\bf x})\Delta t_{n,j}+b(t_n,{\bf x})\Delta t_{n,j}\sqrt{2\Delta t_{n,j}}{\bf a_j}\right){\bf a_j}+R^{GH}_L(YW),
\end{equation*}
where $\E^{{\bf x},h}_{t_n}\left[\cdot\right]$ denotes the approximation of $\E^{{\bf x}}_{t_n}\left[\cdot\right].$
Finally, by inserting these approximations into \eqref{eq:FSchemeYI} and \eqref{eq:FSchemeZI} we obtain
\begin{equation}\label{eq:FSchemeYIF}
\begin{split}
Z^n=(\alpha_{k,0}&+\alpha_{k,1})\tbi{\E}{t_n}{x,h}\left[LI^{n+1}_{h,X^{n,j}}\tbi{Y}{}{n+1}\Delta W_{n,1}^{\top}\right]\\
&+(\alpha_{k,0}+\alpha_{k,1}+\alpha_{k,2})\tbi{\E}{t_n}{x,h}\left[LI^{n+2}_{h,X^{n,j}}\tbi{Y}{}{n+2}\Delta W_{n,2}^{\top}\right]\\
&+\displaystyle\sum_{j=3}^{k}\left(\alpha_{k,j-3}+\alpha_{k,j-2}+\alpha_{k,j-1}+\alpha_{k,j}\right)\tbi{\E}{t_n}{x,h}\left[LI^{n+j}_{h,X^{n,j}}\tbi{Y}{}{n+j}\Delta W_{n,j}^{\top}\right]
\\&+(\alpha_{k,k-2}+\alpha_{k,k-1}+\alpha_{k,k})\tbi{\E}{t_n}{x,h}\left[\tbi{LI^{n+k+1}_{h,X^{n,j}}Y}{}{n+k+1}\Delta W_{n,k+1}^{\top}\right]\\
&+(\alpha_{k,k-1}+\alpha_{k,k})\tbi{\E}{t_n}{x,h}\left[LI^{n+k+2}_{h,X^{n,j}}\tbi{Y}{}{n+k+2}\Delta W_{n,k+2}^{\top}\right]\\
&+\alpha_{k,k}\tbi{\E}{t_n}{x,h}\left[\tbi{LI^{n+k+3}_{h,X^{n,j}}Y}{}{n+k+3}\Delta W_{n,k+3}^{\top}\right]+R_{z,n}^{k,LI_h}+R_{z,n}^{k,\E},
\end{split}
\end{equation}
\begin{equation}\label{eq:FSchemeZIF}
\begin{split}
\alpha_{k,0}Y^n=&-(\alpha_{k,0}+\alpha_{k,1})\tbi{\E}{t_n}{x,h}\left[LI^{n+1}_{h,X^{n,j}}Y^{n+1}\right]-(\alpha_{k,0}+\alpha_{k,1}+\alpha_{k,2})\tbi{\E}{t_n}{x,h}\left[LI^{n+2}_{h,X^{n,j}}Y^{n+2}\right]\\
&-\displaystyle\sum_{j=3}^{k}(\alpha_{k,j-3}+\alpha_{k,j-2}+\alpha_{k,j-1}+\alpha_{k,j})\tbi{\E}{t_n}{x,h}\left[LI^{n+j}_{h,X^{n,j}}Y^{n+j}\right]\\
&-(\alpha_{k,k-2}+\alpha_{k,k-1}+\alpha_{k,k})\tbi{\E}{t_n}{x,h}\left[LI^{n+k+1}_{h,X^{n,j}}Y^{n+k+1}\right]\\
&-(\alpha_{k,k-1}+\alpha_{k,k})\tbi{\E}{t_n}{x,h}\left[LI^{n+k+2}_{h,X^{n,j}}Y^{n+k+2}\right]\\
&-\alpha_{k,k}\tbi{\E}{t_n}{x,h}\left[LI^{n+k+3}_{h,X^{n,j}}Y^{n+k+3}\right]-f(t_n,x,Y^n,Z^n)+R_{y,n}^{k,LI_h}+R_{y,n}^{k,\E}.
\end{split}	
\end{equation}
\begin{remark}
1. The estimate of $R_{y,n}^{k,\E}$ or $R_{z,n}^{k,\E}$ reads \cite{Abramowitz1972, Shen2011, Zhao2014a}
\begin{equation*}
\mathcal{O}\left(\frac{L!}{2^L(2L)!}\right).
\end{equation*}
2. For the local interpolation errors $R_{y,n}^{k,LI_h}$ or $R_{z,n}^{k,LI_h}$ the following estimate holds
\begin{equation}\label{eq:LocalEstimates}
	\mathcal{O}\left(h^{r+1}\right)
\end{equation}
when using $r$-degree polynomial interpolation in $k$-step scheme, and provided that $a,$ $b,$ $f$ and $g$ are sufficiently smooth such that
$\mathcal{L}^{k+4}_{t,x}u(t,x)$  is bounded and $u(t,\cdot) \in C_b^{r+1},$ see \cite{Abramowitz1972, Burden2001,Butcher2008, Zhao2014a}.\\
3. To balance the time discretization error $R^k_{y,n}=\mathcal{O}(\Delta t)^k$ and $R^k_{z,n}=\mathcal{O}(\Delta t)^k,$
one needs to control well both the interpolation and integration error mentioned in last two points.\\
4. For a $k$-step scheme we need to know the support values of $Y^{N_T-i}$ and $Z^{N_T-i}, i=0,\cdots,k+2.$
One can use the following three ways to deal with this problem:
before running the multi-step scheme, we choose a quite smaller $\Delta t$ and run one-step scheme until $N_T-k-2;$
Alternatively, one can prepare these initial values ``iteratively'', namely we compute $Y^{N_T-1}$ and $Z^{N_T-1}$
based on $Y^{N_T}$ and $Z^{N_T}$ with $k=1,$ and the compute $Y^{N_T-2}$ and $Z^{N_T-2}$
based on  $Y^{N_T}, Y^{N_T-1}, Z^{N_T}, Z^{N_T-1}$ with $k=2$ and so on; Finally, one can use the Runge-Kutta scheme proposed in \cite{Crisan2014} with small $\Delta t$ to initialize our proposed multi-step scheme.
\end{remark}
By removing all the error terms, from \eqref{eq:FSchemeYIF} and \eqref{eq:FSchemeZIF} we obtain our fully discrete scheme as follows.
\begin{scheme}
	Assume that $Y^{N_T-i}$ and $Z^{N_T-i}$ on $\mathcal{R}^{N_T-i}_h$ are known for $i=0,1,\cdots,k+2.$ For $n=N_T-k-3,\cdots, 0$ and $x\in \mathcal{R}_h^n,$ $Y^{n}=Y^{n}(x)$ and $Z^{n}=Z^{n}(x)$ can be solved by
	\begin{equation*}
	\hspace{-2.2cm}X^{n,j}=x + a(t_n,x)\Delta t_{n,j}+ b(t_n,x)\Delta W_{n,j},\quad j=1,\cdots,k+3,
	\end{equation*}
	\begin{equation*}
	\begin{split}
	Z^n=(\alpha_{k,0}&+\alpha_{k,1})\tbi{\E}{t_n}{x,h}\left[LI^{n+1}_{h,X^{n,j}}\tbi{Y}{}{n+1}\Delta W_{n,1}^{\top}\right]\\
	&+(\alpha_{k,0}+\alpha_{k,1}+\alpha_{k,2})\tbi{\E}{t_n}{x,h}\left[LI^{n+2}_{h,X^{n,j}}\tbi{Y}{}{n+2}\Delta W_{n,2}^{\top}\right]\\
	&+\displaystyle\sum_{j=3}^{k}\left(\alpha_{k,j-3}+\alpha_{k,j-2}+\alpha_{k,j-1}+\alpha_{k,j}\right)\tbi{\E}{t_n}{x,h}\left[LI^{n+j}_{h,X^{n,j}}\tbi{Y}{}{n+j}\Delta W_{n,j}^{\top}\right]
	\\&+(\alpha_{k,k-2}+\alpha_{k,k-1}+\alpha_{k,k})\tbi{\E}{t_n}{x,h}\left[\tbi{LI^{n+k+1}_{h,X^{n,j}}Y}{}{n+k+1}\Delta W_{n,k+1}^{\top}\right]\\
	&+(\alpha_{k,k-1}+\alpha_{k,k})\tbi{\E}{t_n}{x,h}\left[LI^{n+k+2}_{h,X^{n,j}}\tbi{Y}{}{n+k+2}\Delta W_{n,k+2}^{\top}\right]\\
	&+\alpha_{k,k}\tbi{\E}{t_n}{x,h}\left[\tbi{LI^{n+k+3}_{h,X^{n,j}}Y}{}{n+k+3}\Delta W_{n,k+3}^{\top}\right],
	\end{split}
	\end{equation*}
	\begin{equation*}
	\begin{split}
	\alpha_{k,0}Y^n=&-(\alpha_{k,0}+\alpha_{k,1})\tbi{\E}{t_n}{x,h}\left[LI^{n+1}_{h,X^{n,j}}Y^{n+1}\right]-(\alpha_{k,0}+\alpha_{k,1}+\alpha_{k,2})\tbi{\E}{t_n}{x,h}\left[LI^{n+2}_{h,X^{n,j}}Y^{n+2}\right]\\
	&-\displaystyle\sum_{j=3}^{k}(\alpha_{k,j-3}+\alpha_{k,j-2}+\alpha_{k,j-1}+\alpha_{k,j})\tbi{\E}{t_n}{x,h}\left[LI^{n+j}_{h,X^{n,j}}Y^{n+j}\right]\\
	&-(\alpha_{k,k-2}+\alpha_{k,k-1}+\alpha_{k,k})\tbi{\E}{t_n}{x,h}\left[LI^{n+k+1}_{h,X^{n,j}}Y^{n+k+1}\right]\\
	&-(\alpha_{k,k-1}+\alpha_{k,k})\tbi{\E}{t_n}{x,h}\left[LI^{n+k+2}_{h,X^{n,j}}Y^{n+k+2}\right]\\
	&-\alpha_{k,k}\tbi{\E}{t_n}{x,h}\left[LI^{n+k+3}_{h,X^{n,j}}Y^{n+k+3}\right]-f(t_n,x,Y^n,Z^n).
	\end{split}	
	\end{equation*}
\end{scheme}	
\section{Numerical schemes for coupled FBSDEs}
The authors in \cite{Zhao2014a} extended their scheme proposed for solving decoupled FBSDEs to the one which can solve fully coupled FBSDEs.
Similarly, our Scheme 2 can be extended to solve \eqref{eq:fbsde} in a fully coupled case.
\begin{scheme}
	Assume that $Y^{N_T-i}$ and $Z^{N_T-i}$ on $\mathcal{R}^{N_T-i}_h$ are known for $i=0,1,\cdots,k+2.$ For $n=N_T-k-3,\cdots, 0$ and $x\in \mathcal{R}_h^n,$ $Y^{n}=Y^{n}(x)$ and $Z^{n}=Z^{n}(x)$ can be solved by\\
1. set $Y^{n,0}=Y^{n+1}(x)$ and $Z^{n,0}=Z^{n+1}(x),$ and set $l=0;$\\
2. for $l=0,1,\cdots,$ solve $Y^{n,l+1}=Y^{n,l+1}(x)$ and $Z^{n,l+1}=Z^{n,l+1}(x)$ by 
	\begin{equation*}
	X^{n,j}=x + a(t_n,x,Y^{n,l}(x),Z^{n,l}(x))\Delta t_{n,j}+ b(t_n,x,Y^{n,l}(x),Z^{n,l}(x))\Delta W_{n,j},\quad j=1,\cdots,k+3,
	\end{equation*}
	\begin{equation*}
	\begin{split}
	Z^n&=(\alpha_{k,0}+\alpha_{k,1})\tbi{\E}{t_n}{x,h}\left[LI^{n+1}_{h,X^{n,j}}\tbi{Y}{}{n+1}\Delta W_{n,1}^{\top}\right]
	+(\alpha_{k,0}+\alpha_{k,1}+\alpha_{k,2})\tbi{\E}{t_n}{x,h}\left[LI^{n+2}_{h,X^{n,j}}\tbi{Y}{}{n+2}\Delta W_{n,2}^{\top}\right]\\
	&+\displaystyle\sum_{j=3}^{k}\left(\alpha_{k,j-3}+\alpha_{k,j-2}+\alpha_{k,j-1}+\alpha_{k,j}\right)\tbi{\E}{t_n}{x,h}\left[LI^{n+j}_{h,X^{n,j}}\tbi{Y}{}{n+j}\Delta W_{n,j}^{\top}\right]
	\\&+(\alpha_{k,k-2}+\alpha_{k,k-1}+\alpha_{k,k})\tbi{\E}{t_n}{x,h}\left[\tbi{LI^{n+k+1}_{h,X^{n,j}}Y}{}{n+k+1}\Delta W_{n,k+1}^{\top}\right]\\
	&+(\alpha_{k,k-1}+\alpha_{k,k})\tbi{\E}{t_n}{x,h}\left[LI^{n+k+2}_{h,X^{n,j}}\tbi{Y}{}{n+k+2}\Delta W_{n,k+2}^{\top}\right]
	+\alpha_{k,k}\tbi{\E}{t_n}{x,h}\left[\tbi{LI^{n+k+3}_{h,X^{n,j}}Y}{}{n+k+3}\Delta W_{n,k+3}^{\top}\right],
	\end{split}
	\end{equation*}
	\begin{equation*}
	\begin{split}
	\alpha_{k,0}Y^n&=-(\alpha_{k,0}+\alpha_{k,1})\tbi{\E}{t_n}{x,h}\left[LI^{n+1}_{h,X^{n,j}}Y^{n+1}\right]-(\alpha_{k,0}+\alpha_{k,1}+\alpha_{k,2})\tbi{\E}{t_n}{x,h}\left[LI^{n+2}_{h,X^{n,j}}Y^{n+2}\right]\\
	&-\displaystyle\sum_{j=3}^{k}(\alpha_{k,j-3}+\alpha_{k,j-2}+\alpha_{k,j-1}+\alpha_{k,j})\tbi{\E}{t_n}{x,h}\left[LI^{n+j}_{h,X^{n,j}}Y^{n+j}\right]\\
	&-(\alpha_{k,k-2}+\alpha_{k,k-1}+\alpha_{k,k})\tbi{\E}{t_n}{x,h}\left[LI^{n+k+1}_{h,X^{n,j}}Y^{n+k+1}\right]
	-(\alpha_{k,k-1}+\alpha_{k,k})\tbi{\E}{t_n}{x,h}\left[LI^{n+k+2}_{h,X^{n,j}}Y^{n+k+2}\right]\\
	&-\alpha_{k,k}\tbi{\E}{t_n}{x,h}\left[LI^{n+k+3}_{h,X^{n,j}}Y^{n+k+3}\right]-f(t_n,x,Y^n,Z^n).
	\end{split}	
	\end{equation*}
	until $\max\left(\left|Y^{n,l+1}-Y^{n,l}\right|,\left|Z^{n,l+1}-Z^{n,l}\right|\right)<\epsilon_0,$\\
3. let $Y^n=Y^{n,l+1}$ and $Z^n=Z^{n,l+1}.$
\end{scheme}	
\begin{remark}
1. Scheme 3 coincides with Scheme 2 if $a$ and $b$ do not depend on $Y$ and $Z.$\\
2. We only assume that the coupled FBSDEs are uniquely solvable, the lacking analysis will be the task of future work.
\end{remark}
\section{Numerical experiments}
In this section we use some numerical examples to show that our Schemes 2 and 3 can reach ninth-order convergence rate for
solving FBSDEs. The uniform partitions in both time and space will be used, that is, the time interval $[0,T]$ will be uniformly
divided into $N_T$ parts with $\Delta t=\frac{T}{N_T}$ such that $t_n=n\Delta t, n=0,1,\cdots,N_T;$ the space partition is $\mathcal{R}_h^n=R_h$ for all $n$ with
\begin{equation*}
	\mathcal{R}_h=\mathcal{R}_{1,h}\times \mathcal{R}_{2,h} \times \cdots \mathcal{R}_{n,h},
\end{equation*}
where $\mathcal{R}_{j,h}$ is the partition of $\mathbb{R}$
\begin{equation*}
\mathcal{R}_{j,h}=\left\{x_i^j: x_i^j=ih, i=0,\pm 1,\cdots,\pm \infty \right\},\quad j=1,2,\cdots,n.
\end{equation*}

In our numerical experiments we choose the local Lagrange interpolation for $LI^n_{h,x}$ based on the set of
some neighbor grids near $x,$ i.e., $\mathcal{R}_{h,x} \subset \mathcal{R}_h$ such that  \eqref{eq:LocalEstimates} holds.
Following \cite{Zhao2014}, we set sufficiently many Gauss-Hermite quadrature points such that the quadrature error could be negligible.
Note that the truncation error is defined in \eqref{eq:errorTerm},
in order to thus balance the time and space truncation error in our numerical examples, we force $h^{r+1}=(\Delta t)^{k+1},$ where $r$
is the degree of the Lagrangian interpolation polynomials. For example, one can firstly specify a value of $r,$ and then adjust the value of
$h$ such that $h=\Delta t^{\frac{k+1}{r+1}}.$ For the numerical results in this paper, $r$ is set to be a value from the set $\{10, 11, \cdots, 21\}$
to control the errors. Furthermore, we will consider $k$ from $3$ such that at least one combination of four $\alpha_{k,i}\,s$ is included, but
until $k=9$ due to the stability condition, see Scheme 2 and 3. Finally, CR and RT are used to denote the convergence rate and the running time
in second, respectively. For the comparison purpose, we will directly take examples considered in \cite{Zhao2014}. 
Numerical experiment were performed in MATLAB with an Intel(R) Core(TM) i5-8350 CPU @ 1.70 GHz and 15 G RAM.

\paragraph{Example 1}
The first example reads
\begin{equation*}
\begin{cases}
&dX_t=\frac{1}{1+2\exp(t+X_t)}\,dt + \frac{\exp(t+X_t)}{1+\exp(t+X_t)}\,dW_t,\quad X_0=1,\\
&-dY_t=\left(-\frac{2Y_t}{1+2\exp(t+X_t)}-\frac{1}{2}\left(\frac{Y_tZ_t}{1+\exp(t+X_t)}-Y_t^2Z_t\right)\right) \,dt-Z_t\,dW_t,\\
&Y_T=\frac{\exp(T+X_T)}{1+\exp(T+X_T)},
\end{cases}
\end{equation*}
with the analytic solution
\begin{equation*}
\begin{cases}
&Y_t=\frac{\exp(t+X_t)}{1+\exp(t+X_t)},\\
&Z_t=\frac{(\exp(t+X_t))^2}{(1+\exp(t+X_t))^3}.
\end{cases}
\end{equation*}
Obviously, in this example, the generator $a$ and $b$ does not depend on $Y_t$ and
$Z_t,$ i.e., a decoupled FBSDE. All the convergence rates, running time and absolute errors are reported
in Table \ref{tab:exp01}. In the proposed scheme, $k+3$ points are needed for the iterations, i.e., one need
$12$ points when $k=9.$ From the other side we can not choose a large value for $N_T$ due to the accuracy of double
precision. Therefore, to show the convergence rate up to ninth order we consider $N_T=\{16, 20, 24, 28, 32 \}$ in this example.
\begin{table}[h!]
\centering
\begin{tabular}{|c|c|c|c|c|c|c|c|}
\hline
Scheme 2 && $N_T=16$ & $N_T=20$ & $N_T=24$ & $N_T=28$ & $N_T=32$ & CR \\
\hline
\multirow{3}{*}{$k=3$}&
$|Y^0-Y_0|$ & 4.717e-06 &2.613e-06& 1.569e-06& 1.015e-06& 6.905e-07 & 2.78 \\ \cline{2-8}

&$|Z^0-Z_0|$ &  2.547e-05&1.552e-05 & 1.009e-05 &   6.889e-06  &  4.903e-06  &  2.39 \\\cline{2-8}

&RT &0.37 &0.49  &  0.65 &   0.91&    1.19  & \\\cline{2-8}

\hline
\multirow{3}{*}{$k=4$}&
$|Y^0-Y_0|$ &6.871e-07& 3.152e-07&1.629e-07 &9.240e-08 &5.618e-08  & 3.61 \\ \cline{2-8}

&$|Z^0-Z_0|$ & 6.879e-06 &3.097e-06&1.595e-06&9.027e-07 &5.488e-07 & 3.65 \\\cline{2-8}

&RT &0.39 &0.61&0.85& 1.12&1.41 &\\\cline{2-8}
\hline
\multirow{3}{*}{$k=5$}&
$|Y^0-Y_0|$ & 5.623e-08&2.077e-08&9.011e-09&4.355e-09 &2.343e-09&    4.59 \\ \cline{2-8}

&$|Z^0-Z_0|$ & 6.522e-07&2.427e-07 & 1.047e-07 &   5.016e-08 &   2.704e-08&    4.60 \\\cline{2-8}

&RT & 0.46 & 0.69 &   0.95&   1.25&    1.58&   \\ \cline{2-8}
\hline
\multirow{3}{*}{$k=6$}&
$|Y^0-Y_0|$ & 3.549e-09 &1.073e-09&3.929e-10 &   1.623e-10&    7.549e-11 &   5.56 \\ \cline{2-8}

&$|Z^0-Z_0|$ & 5.658e-08& 1.632e-08&  6.000e-09 &   2.519e-09&    1.168e-09&    5.59 \\\cline{2-8}

&RT &0.52&    0.85 &    1.26 &   1.74  &  2.10&          \\\cline{2-8}
\hline
\multirow{3}{*}{$k=7$}&
$|Y^0-Y_0|$ &2.156e-10& 4.809e-11& 1.457e-11 & 5.075e-12& 2.019e-12 & 6.73\\ \cline{2-8}

&$|Z^0-Z_0|$ &6.349e-09& 1.556e-09& 4.796e-10& 1.749e-10&  7.147e-11&  6.47 \\\cline{2-8}

&RT & 0.60&  1.10 &   1.65&    2.15 &   2.79  \\\cline{2-8}
\hline
\multirow{3}{*}{$k=8$}&
$|Y^0-Y_0|$ &6.025e-11 & 8.573e-12&3.292e-12& 7.027e-13& 4.868e-13 &  7.10 \\ \cline{2-8}

&$|Z^0-Z_0|$ &  1.029e-09& 1.934e-10 &  5.811e-11&    1.459e-11  &6.696e-12&    7.35\\\cline{2-8}

&RT & 0.62 &   1.17 &   1.73&    2.37&    3.11 &   \\\cline{2-8}
\hline
\multirow{3}{*}{$k=9$}&
$|Y^0-Y_0|$ &  2.315e-11&    4.131e-12&    9.073e-13&    2.169e-13 &   2.398e-14 &   9.55\\ \cline{2-8}

&$|Z^0-Z_0|$ & 3.672e-10&    5.073e-11&    1.184e-11&    2.528e-12 &   5.760e-13 &   9.19 \\\cline{2-8}

&RT &0.69  &1.32 &   2.04&    2.83&    3.74& \\\cline{2-8}
\hline
\end{tabular}
\caption{Errors, running time and convergence rates for Example 1, $T=1$}\label{tab:exp01}
\end{table}
From Table \ref{tab:exp01} we see that the quite high accuracy of Scheme 2 for solving decoupled
FBSDEs. Scheme 2 is a $k$-order scheme up to $k=9,$ and more efficient for taking a larger value for $k,$
which is consistent with the theory \cite{Butcher2008}, see also Table \ref{tab:4}

For the second example we consider the coupled FBSDE (taken from \cite{Zhao2014}) to test Scheme 3, in which
an iterative process is required with longer computational time. 
\paragraph{Example 2}
 \begin{equation*}
\begin{cases}
&dX_t=-\frac{1}{2}\sin(t+X_t)\cos(t+X_t)(Y^2_t+Z_t)\,ds\\
&\qquad+\frac{1}{2}\cos(t+X_t)+(Y_t\sin(t+X_t)+Z_t+1)\,dW_s,\quad X_0=1.5,\\
&-dY_t=\left(Y_t Z_t-\cos(t+X_t)\right)\,dt-Z_t\,dW_t,\\
&Y_T=\sin(T+X_T).
\end{cases}
\end{equation*}
has the analytic solution
\begin{equation*}
\begin{cases}
&Y_t=\sin(t+X_t),\\
&Z_t=\cos^2(t+X_t).
\end{cases}
\end{equation*}
In this coupled FBSDE, the diffusion coefficient $b$ depends on $X, Y$ and $Z,$ i.e., quite general.
Due to the same reasons as those explained for Example 1, we set $N_T=\{13, 15, 17, 19, 21\}$ in order to show the convergence rate up to
ninth order.
\begin{table}[h!]
	\centering
	\begin{tabular}{|c|c|c|c|c|c|c|c|}
		\hline
		Scheme 3 && $N_T=13$ & $N_T=15$ & $N_T=17$ & $N_T=19$ & $N_T=21$ & CR \\
		\hline
		\multirow{3}{*}{$k=3$}&
		$|Y^0-Y_0|$ & 2.269e-04 & 1.398e-04 & 9.186e-05 &   6.025e-05 &   4.336e-05 &   3.47 \\ \cline{2-8}
		
		&$|Z^0-Z_0|$ &  1.562e-04    &1.143e-04   & 6.555e-05  &  4.431e-05  &  3.359e-05  &  3.36\\\cline{2-8}
		
		&RT &  3.89& 4.79&  5.78 &  6.53 &   7.84   &    \\\cline{2-8}
		
		\hline
		\multirow{3}{*}{$k=4$}&
		$|Y^0-Y_0|$ & 9.569e-06 &   7.654e-06  &  4.799e-06 &   2.734e-06   & 1.571e-06  &  3.83 \\ \cline{2-8}
		
		&$|Z^0-Z_0|$ &1.447e-04  &  1.126e-04  &  6.268e-05  &  3.009e-05  &  1.223e-05   & 5.14 \\\cline{2-8}
		
		&RT &4.51&    5.73  & 6.79 &   8.12  &  10.43&   \\\cline{2-8}
		\hline
		\multirow{3}{*}{$k=5$}&
		$|Y^0-Y_0|$ & 4.773e-07 &   1.740e-07  &  3.835e-08 &   6.464e-08   & 2.325e-08    &5.96 \\ \cline{2-8}
		
		&$|Z^0-Z_0|$ &2.433e-06  &  6.129e-07  &  2.215e-08  &  1.988e-07    &2.968e-07    &4.89\\\cline{2-8}
		
		&RT & 21.07&    28.68 &   35.43 &    41.14 &  53.55 &     \\\cline{2-8}
		\hline
		\multirow{3}{*}{$k=6$}&
		$|Y^0-Y_0|$ &4.469e-08 &   2.509e-08   & 1.361e-08  &  6.572e-09  &  3.257e-09 &   5.47 \\ \cline{2-8}
		
		&$|Z^0-Z_0|$ &4.257e-07 &   3.145e-07  &  1.629e-07 &   7.017e-08  &  2.121e-08  &  6.15 \\\cline{2-8}
		
		&RT & 33.49&    48.76&   63.83 &    80.93 &    99.07    \\\cline{2-8}
		\hline
		\multirow{3}{*}{$k=7$}&
		$|Y^0-Y_0|$ & 6.510e-10 &   4.904e-10    &1.536e-11  &  4.256e-11  &  3.250e-11 &   7.20\\ \cline{2-8}
		
		&$|Z^0-Z_0|$ &1.218e-08  &  1.207e-08  &  8.072e-10   & 3.704e-09   & 2.249e-10&    7.56 \\\cline{2-8}
		
		&RT &  34.39 &    53.76&    75.34 &    97.09 &   120.69     &      \\\cline{2-8}
		\hline
		\multirow{3}{*}{$k=8$}&
		$|Y^0-Y_0|$ & 1.876e-10 &8.477e-11& 3.098e-11& 1.028e-11&  5.961e-12 &    7.51 \\ \cline{2-8}
		
		&$|Z^0-Z_0|$ & 8.841e-09&2.877e-09& 1.727e-10& 7.179e-11&    5.141e-10 &   8.19 \\\cline{2-8}
		
		&RT & 35.76&    90.17&    149.87 &  172.70   &247.68    \\\cline{2-8}
		\hline
		\multirow{3}{*}{$k=9$}&
		$|Y^0-Y_0|$ &  1.926e-11 & 5.744e-12& 5.828e-13& 9.910e-13&  3.098e-14&    12.10 \\ \cline{2-8}
		
		&$|Z^0-Z_0|$ & 2.502e-10& 1.526e-10&  3.190e-11& 3.039e-11&  2.234e-12&   9.07 \\\cline{2-8}
		
		&RT & 55.19  &  173.80  &  275.15&    384.02&    503.85 &   \\\cline{2-8}
		\hline
	\end{tabular}
	\caption{Errors, running time and convergence rates for Example 2, $T=1$}\label{tab:exp02}
\end{table}
From the results listed in Table \ref{tab:exp02} one can clearly draw same conclusions as those having been for Example 1. 

Finally, we illustrate the accuracy of the proposed scheme for a two-dimensional example, which is also taken from
\cite{Zhao2014} and reads
\paragraph{Example 3}
\begin{equation*}
\left\{
\begin{array}{l}
\begin{pmatrix}
dX^1_t\\
dX^2_t
\end{pmatrix}=\begin{pmatrix}
\frac{1}{2}\sin^2(t+X^1_t)\\
\frac{1}{2}\sin^2(t+X^2_t)
\end{pmatrix}\,dt + \begin{pmatrix}
\frac{1}{2}\cos^2(t+X^1_t)\\
\frac{1}{2}\cos^2(t+X^2_t)
\end{pmatrix}\,dW_t,\quad
\begin{pmatrix}
X^1_0\\
X^2_0
\end{pmatrix}=\begin{pmatrix}
0\\
0
\end{pmatrix},\\
\begin{pmatrix}
dY^1_t\\
dY^2_t
\end{pmatrix}=\begin{pmatrix}
&-\frac{3}{2}\cos(t+X^1_t)\sin(t+X^2_t)-\frac{3}{2}\sin(t+X^1_t)\cos(t+X^2_t)-Z^2_t\\
&+\frac{1}{2}Y_t^1\left(\frac{1}{4}\cos^4(t+X^2_t)+\frac{1}{4}\cos^4(t+X^1_t)\right)-\frac{1}{4}(Y^2_t)^3\\
&\frac{3}{2}\sin(t+X^1_t)\cos(t+X^2_t)+\frac{3}{2}\cos(t+X^1_t)\sin(t+X^2_t)-Z^1_t\\
&+\frac{1}{2}Y_t^2	\left(\frac{1}{4}\cos^4(t+X^2_t)+\frac{1}{4}\cos^4(t+X^1_t)\right)-\frac{1}{4}Y^1_t(Y^2_t)^2
\end{pmatrix}\,dt - \begin{pmatrix}
Z^1_t\\
Z^2_t
\end{pmatrix}\,dW_t,\\
 \begin{pmatrix}
 Y^1_T\\
 Y^2_T
 \end{pmatrix}=\begin{pmatrix}
 \sin(T+X^1_T)\sin(T+X^2_T)\\
 \cos(T+X^1_T)\cos(T+X^2_T)
 \end{pmatrix}
 \end{array}\right.
 \end{equation*}
with the analytic solution
 \begin{equation*}
 \left\{
 \begin{array}{l}
 \begin{pmatrix}
 Y^1_t\\
 Y^2_t
 \end{pmatrix}=\begin{pmatrix}
 \sin(t+X^1_t)\sin(t+X^2_t)\\
 \cos(t+X^1_t)\cos(t+X^2_t)
 \end{pmatrix},  \\
 \begin{pmatrix}
 Z_t^1\\
 Z_t^2
 \end{pmatrix}=\begin{pmatrix}
 &\frac{1}{2}\cos(t+X^1_t)\sin(t+X^2_t)\cos^2(t+X^2_t)\\
 &\qquad\qquad+\frac{1}{2}\sin(t+X^1_t)\cos(t+X^2_t)\cos^2(t+X^1_t)\\
 &-\frac{1}{2}\sin(t+X^1_t)\cos^3(t+X^2_t)-\frac{1}{2}\cos^3(t+X^1_t)\sin(t+X^2_t)
 \end{pmatrix}.
 \end{array}\right.
 \end{equation*}
The numerical approximations are reported in Table \ref{tab:exp03}, which show that our multi-step scheme
is still quite highly accurate for solving a two-dimensional FBSDE.
\begin{table}[h!]
	\centering
	\begin{tabular}{|c|c|c|c|c|c|c|c|}
		\hline
		Scheme 2 && $N_T=16$ & $N_T=20$ & $N_T=24$ & $N_T=28$ & $N_T=32$ & CR \\
		\hline
		\multirow{5}{*}{$k=3$}&
		$|Y^{0,1}-Y^1_0|$ &4.308e-03 & 2.447e-03 &   1.495e-03 &   9.726e-04  &  6.651e-04 &   2.70 \\ \cline{2-8}
		
		&$|Y^{0,2}-Y^2_0|$ & 3.896e-03 &   2.200e-03  &  1.340e-03  &  8.705e-04 &   5.949e-04   & 2.71 \\\cline{2-8}
		
		&$|Z^{0,1}-Z^1_0|$ &3.887e-03 &   2.041e-03  &  1.201e-03& 7.602e-04 &   5.110e-04  &  2.93 \\\cline{2-8}
		&$|Z^{0,2}-Z^2_0|$& 2.980e-03  &  1.491e-03  &  8.465e-04 &   5.243e-04  &  3.451e-04   & 3.11 \\\cline{2-8}
		&RT &  34.65&    53.68&    74.87&    92.42  &  130.14        \\\cline{2-8}
		\hline
	    \multirow{5}{*}{$k=4$}
    	&$|Y^{0,1}-Y^1_0|$ & 8.829e-04  &  4.437e-04 &   2.417e-04&    1.413e-04 &   8.759e-05 &   3.34 \\ \cline{2-8}
	
   	&$|Y^{0,2}-Y^2_0|$ & 7.776e-04 &   3.845e-04   & 2.072e-04  &  1.202e-04  &  7.413e-05 &   3.39 \\\cline{2-8}
	
	&$|Z^{0,1}-Z^1_0|$ &  6.653e-04  &  2.556e-04 &   1.161e-04  &  5.899e-05  &  3.270e-05  &  4.35 \\\cline{2-8}
	&$|Z^{0,2}-Z^2_0|$ &6.845e-04 &   2.659e-04  &  1.218e-04 &   6.268e-05 &   3.521e-05 &   4.29\\\cline{2-8}
	&RT &   43.52 &    70.30&    100.41&    142.95&    195.94        \\\cline{2-8}
		\hline
		\multirow{5}{*}{$k=5$}&
	$|Y^{0,1}-Y^1_0|$ & 5.686e-05 &   2.367e-05  &  1.084e-05   & 5.445e-06   & 2.952e-06 &   4.27 \\ \cline{2-8}
	
	&$|Y^{0,2}-Y^2_0|$ & 5.131e-05 &   2.112e-05  &  9.605e-06 &   4.801e-06 &   2.595e-06  &  4.31 \\\cline{2-8}
	
	&$|Z^{0,1}-Z^1_0|$ & 7.703e-05 &   2.621e-05 &   1.067e-05 &   4.923e-06 &   2.517e-06 &   4.94 \\\cline{2-8}
	&$|Z^{0,2}-Z^2_0|$&7.237e-05 &   2.402e-05 &   9.582e-06 &   4.363e-06   & 2.201e-06 &   5.04\\\cline{2-8}
	&RT &  50.72&    79.97&    132.00 &    186.09&    231.84 &          \\\cline{2-8}
	\hline
		\multirow{5}{*}{$k=6$}&
$|Y^{0,1}-Y^1_0|$ &   4.277e-06& 1.706e-06& 7.258e-07&  3.354e-07&  1.674e-07  &   4.68\\ \cline{2-8}

&$|Y^{0,2}-Y^2_0|$ & 3.857e-06&  1.509e-06&    6.336e-07&    2.901e-07&    1.437e-07&     4.75\\\cline{2-8}

&$|Z^{0,1}-Z^1_0|$ & 5.569e-06&    1.667e-06&    5.803e-07&    2.309e-07&    1.016e-07 &    5.78\\\cline{2-8}
&$|Z^{0,2}-Z^2_0|$&5.531e-06&    1.655e-06&    5.791e-07 &   2.313e-07  &  1.026e-07 &    5.76\\\cline{2-8}
&RT & 64.97 &   110.97&    173.50&    246.96&    337.71&  \\\cline{2-8}
		\hline
			\multirow{5}{*}{$k=7$}&
	$|Y^{0,1}-Y^1_0|$ &5.753e-07 &1.843e-07& 6.513e-08&  2.568e-08&    1.115e-08&     5.70 \\ \cline{2-8}
	
	&$|Y^{0,2}-Y^2_0|$ &  5.312e-07&    1.669e-07 &   5.836e-08 &   2.283e-08  &  9.864e-09&     5.76\\\cline{2-8}
	
	&$|Z^{0,1}-Z^1_0|$ &9.466e-07 &   2.291e-07 &   7.027e-08 &   2.528e-08 &   1.027e-08&     6.52\\\cline{2-8}
	&$|Z^{0,2}-Z^2_0|$& 9.091e-07 &   2.157e-07 &   6.529e-08 &   2.326e-08 &   9.404e-09 &    6.59 \\\cline{2-8}
	&RT &  73.89&    139.46&    234.34&    339.44&    453.70&        \\\cline{2-8}
		\hline
			\multirow{5}{*}{$k=8$}&
	$|Y^{0,1}-Y^1_0|$ & 2.586e-08 &   1.080e-08  &  4.077e-09 &   1.592e-09  &  6.652e-10   &  5.30\\ \cline{2-8}
	
	&$|Y^{0,2}-Y^2_0|$ &  2.384e-08&    9.742e-09&    3.623e-09 &   1.400e-09   & 5.802e-10&     5.39 \\\cline{2-8}
	
	&$|Z^{0,1}-Z^1_0|$ &6.238e-08 &   1.548e-08  &  4.186e-09 &   1.299e-09  &  4.791e-10 &    7.06 \\\cline{2-8}
	&$|Z^{0,2}-Z^2_0|$& 6.290e-08  &  1.534e-08  &  4.077e-09 &   1.255e-09  &  4.579e-10 &    7.14 \\\cline{2-8}
	&RT &  84.65&   177.96&    299.48&   436.73&    592.02&   \\\cline{2-8}
		\hline
		\multirow{5}{*}{$k=9$}&
$|Y^{0,1}-Y^1_0|$ & 8.294e-09&    2.371e-09&    6.964e-10  &  2.250e-10  &  8.045e-11&     6.70 \\ \cline{2-8}

&$|Y^{0,2}-Y^2_0|$ &7.739e-09&    2.153e-09&    6.229e-10 &   1.991e-10  &  6.989e-11&     6.80 \\\cline{2-8}

&$|Z^{0,1}-Z^1_0|$ & 2.771e-08&    5.329e-09 &   1.262e-09&    3.491e-10 &   1.245e-10 &    7.84\\\cline{2-8}
&$|Z^{0,2}-Z^2_0|$&2.748e-08  &  5.202e-09   & 1.217e-09  &  3.350e-10  &  9.998e-11  &   8.08 \\\cline{2-8}
&RT & 91.19 &   210.99 &    356.54&    534.76&    760.94&      \\\cline{2-8}
		\hline
	\end{tabular}
	\caption{Errors, running time and convergence rates for Example 3, $T=1$}\label{tab:exp03}
\end{table}
We observe that the convergence rates are roughly consistent with the theoretical results,
the slight deviation comes from the quadratures and especially the two-dimensional interpolations.
Obviously, the high efficiency and accuracy have been shown in this example. Note that the parallel
computing toolbox in MATLAB has been used in this example, more precisely, the parallel for-Loops (parfor)
is used for the two-dimensional interpolation on the grid points.

\section{Conclusion}
In this work, by combining the multi-steps we have adopted the high-order multi-step method in
[W. Zhao, Y. Fu and T. Zhou, SIAM J. Sci. Comput., 36(4) (2014), pp.A1731-A1751] for solving FBSDEs.
First of all, our new schemes allow for higher convergence rate up to ninth order, and are more efficient.
Secondly, they keep the key feature of the method in [W. Zhao, Y. Fu and T. Zhou, SIAM J. Sci. Comput. 36(4), pp.A1731-A1751],
that is the numerical solution of backward component maintains the higher-order accuracy by using the Euler method to the forward component.  This feature makes our schemes be promising in solving problems in practice.
The effectiveness and higher-order accuracy have been confirmed by the numerical experiments.
A rigorous stability analysis for the proposed schemes is the task of future work.

\section*{Acknowledge}
We thank Suman Kumar from the University of Wuppertal for his assistance with Matlab programming partially used for the numerical experiments.
\bibliography{mybibfile}

\begin{thebibliography}{}

\bibitem[Abramowitz and Stegun, 1972]{Abramowitz1972}
Abramowitz, M. and Stegun, I. (1972).
\newblock {\em Handbook of Mathematical Functions}.
\newblock Dover Publications.
\newblock Dover Books on Mathematics.

\bibitem[Bender and Steiner, 2012]{Bender2012}
Bender, C. and Steiner, J. (2012).
\newblock Least-squares monte carlo for backward sdes.
\newblock {\em Numer. Methods Finance}, 12:257--289.

\bibitem[Bender and Zhang, 2008]{Bender2008}
Bender, C. and Zhang, J. (2008).
\newblock Time discretization and markovian iteration for coupled fbsdes.
\newblock {\em Ann. Appl. Probab.}, 18:143--177.

\bibitem[Bouchard and Touzi, 2004]{Bouchard2004}
Bouchard, B. and Touzi, N. (2004).
\newblock Discrete-time approximation and monte-carlo simulation of backward
  stochastic differential equations.
\newblock {\em Stoch. Proc. Appl.}, 111:175--206.

\bibitem[Burden and Faires, 2001]{Burden2001}
Burden, R.~L. and Faires, J.~D. (2001).
\newblock {\em Numerical Analysis}.
\newblock Higher Education Press/Cengage Learning.
\newblock 7th ed.

\bibitem[Butcher, 2008]{Butcher2008}
Butcher, J.~C. (2008).
\newblock {\em Numerical methods for ordinary differential equations}.
\newblock John Wiley, Chichester, UK.

\bibitem[Crisan and Chassagneux, 2014]{Crisan2014}
Crisan, D. and Chassagneux, J.~F. (2014).
\newblock Runge-kutta schemes for backward stochastic differential equations.
\newblock {\em Ann. Appl. Probab.}, 24:679--720.

\bibitem[Crisan and Manolarakis, 2010]{Crisan2010}
Crisan, D. and Manolarakis, K. (2010).
\newblock Solving backward stochastic differential equations using the cubature
  method: Application to nonlinear pricing.
\newblock {\em SIAM J. FINAN. MATH}, 3(1):534--571.

\bibitem[Cvitanic and Zhang, 2006]{Cvitanic2006}
Cvitanic, J. and Zhang, J. (2006).
\newblock The steepest descent method for forward-backward sdes.
\newblock {\em Electron. J. Probab.}, 16:940--968.

\bibitem[Delarue and Menozzi, 2006]{Delarue2006}
Delarue, F. and Menozzi, S. (2006).
\newblock A forward-backward stochastic algorithm for quasi-linear pdes.
\newblock {\em Ann. Appl. Probab.}, 16(1):140--184.

\bibitem[Douglas et~al., 1996]{Douglas1996}
Douglas, J., Ma, J., and Protter, P. (1996).
\newblock Numerical methods for forward-backward stochastic differential
  equations.
\newblock {\em Ann. Appl. Probab.}, 6:940--968.

\bibitem[Fornberg, 1988]{Fornberg1988}
Fornberg, B. (1988).
\newblock Generation of finite difference formulas on arbitrarily spaced grids.
\newblock {\em Math. Comput.}, 51(184):699--706.

\bibitem[Fu et~al., 2017]{Fu2017}
Fu, Y., Zhao, W., and Zhou, T. (2017).
\newblock Efficient spectral sparse grid approximations for solving
  multi-dimensional forward backward sdes.
\newblock {\em Discrete Cont. Dyn-B.}, 22(9):3439--3458.

\bibitem[Gobet et~al., 2005]{Gobet2005}
Gobet, E., Lemor, J.~P., and Warin, X. (2005).
\newblock A regression-based monte carlo method to solve backward stochastic
  differential equations.
\newblock {\em Ann. Appl. Probab.}, 15:2172--2202.

\bibitem[Lemor et~al., 2006]{Lemor2006}
Lemor, J., Gobet, E., and Warin, X. (2006).
\newblock Rate of convergence of an empirical regression method for solving
  generalized backward stochastic differential equations.
\newblock {\em Bernoulli}, 12:889--916.

\bibitem[Lepeltier and Martin, 1997]{Lepeltier1997}
Lepeltier, J.~P. and Martin, J.~S. (1997).
\newblock Backward stochastic differential equations with continuous generator.
\newblock {\em Statist. Probab. Lett.}, 32(425--430).

\bibitem[Ma et~al., 1994]{Ma1994}
Ma, J., Protter, P., and Yong, J. (1994).
\newblock Solving forward-backward stochastic differential equations
  explicity-a four step scheme.
\newblock {\em Probab. Theory Related Fields}, 98(3):339--359.

\bibitem[Ma et~al., 2008]{Ma2008}
Ma, J., Shen, J., and Zhao, Y. (2008).
\newblock On numerical approximations of forward-backward stochastic
  differential equations.
\newblock {\em SIAM J. Numer. Anal.}, 46(5):2636--2661.

\bibitem[Ma and Zhang, 2005]{Ma2005}
Ma, J. and Zhang, J. (2005).
\newblock Representations and regularities for solutions to bsdes with
  reflections.
\newblock {\em Stoch. Proc. Appl.}, 115:539--569.

\bibitem[Milsetin and Tretyakov, 2006]{Milstein2006}
Milsetin, G.~N. and Tretyakov, M.~V. (2006).
\newblock Numerical algorithms for forward-backward stochastic differential
  equations.
\newblock {\em SIAM J. SCI. COMPUT.}, 28:561--582.

\bibitem[Pardoux and Peng, 1990]{Pardoux1990}
Pardoux, E. and Peng, S. (1990).
\newblock Adapted solution of a backward stochastic differential equations.
\newblock {\em System and Control Letters}, 14:55--61.

\bibitem[Pardoux and Peng, 1992]{Pardoux1992}
Pardoux, E. and Peng, S. (1992).
\newblock Backward stochastic differential equation and quasilinear parabolic
  partial differential equations.
\newblock {\em Lectures Notes in CSI.}, 176:200--217.

\bibitem[Peng, 1991]{Peng1991}
Peng, S. (1991).
\newblock Probabilistic interpretation for systems of quasilinear parabolic
  partial differential equations.
\newblock {\em Stochastics and Stochastic Reports}, 37(1--2):61--74.

\bibitem[Peng and Wu, 1999]{Peng1999}
Peng, S. and Wu, Z. (1999).
\newblock Fully coupled forward-backward stochastic differential equations and
  applications to optimal control.
\newblock {\em SIAM J. Control Optim.}, 37:825--843.

\bibitem[Ruijter and Oosterlee, 2015]{Ruijter2015}
Ruijter, M.~J. and Oosterlee, C.~W. (2015).
\newblock A fourier cosine method for an efficient computation of solutions to
  bsdes.
\newblock {\em SIAM J. SCI. COMPUT.}, 37(2):A859--A889.

\bibitem[Shen et~al., 2011]{Shen2011}
Shen, J., Tang, T., and Wang, L. (2011).
\newblock {\em Spectral Methods: Algorithms, Analysis and Applications}.
\newblock Springer-Verlag, Berlin.

\bibitem[Teng, 2019]{Teng2019}
Teng, L. (2019).
\newblock A review of tree-based approaches to solve forward-backward
  stochastic differential equations.
\newblock {\em arXiv:1809.00325v4, available on webpage at
  \url{https://arxiv.org/pdf/1809.00325v4.pdf}}.

\bibitem[Teng et~al., 2020]{Teng2020}
Teng, L., Lapitckii, A., and G{\"unther}, M. (2020).
\newblock A multi-step scheme based on cubic spline for solving backward
  stochastic differential equations.
\newblock {\em Appl. Numer. Math.}, 150.

\bibitem[Zhang et~al., 2013]{Zhang2013}
Zhang, G., Gunzburger, M., and Zhao, W. (2013).
\newblock A sparse-grid method for multi-dimensional backward stochastic
  differential equations.
\newblock {\em J. Comput. Math.}, 31(3):221--248.

\bibitem[Zhang, 2001]{Zhang2001}
Zhang, J. (2001).
\newblock {\em Some fine properties of backward stochastic differential
  equations}.
\newblock PhD thesis, Purdue University, West Lafayette, IN.

\bibitem[Zhang, 2004]{Zhang2004}
Zhang, J. (2004).
\newblock A numerical scheme for bsdes.
\newblock {\em Ann. Appl. Probab.}, 14:459--488.

\bibitem[Zhao et~al., 2006]{Zhao2006}
Zhao, W., Chen, L., and Peng, S. (2006).
\newblock A new kind of accurate numerical method for backward stochastic
  differential equations.
\newblock {\em SIAM J. SCI. COMPUT.}, 28(4):1563--1581.

\bibitem[Zhao et~al., 2014a]{Zhao2014}
Zhao, W., Fu, Y., and Zhou, T. (2014a).
\newblock New kinds of high-order multistep schemes for coupled forward
  backward stochastic differential equations.
\newblock {\em SIAM J. SCI. COMPUT.}, 36(4):A1731--A1751.

\bibitem[Zhao et~al., 2013]{Zhao2013}
Zhao, W., Li, Y., and Ju, L. (2013).
\newblock Error estimates of the crank-nicolson scheme for solving backward
  stochastic differential equations.
\newblock {\em Int. J. Numer. Anal. Mode.l}, 10(4):876--898.

\bibitem[Zhao et~al., 2012]{Zhao2012}
Zhao, W., Li, Y., and Zhang, G. (2012).
\newblock A generalized $\theta$-scheme for solving backward stochastic
  differential equations.
\newblock {\em Discrete Cont. Dyn-B.}, 17(5):1585--1603.

\bibitem[Zhao et~al., 2009]{Zhao2009}
Zhao, W., Wang, J., and Peng, S. (2009).
\newblock Error estimates of the theta-scheme for backward stochastic
  differential equations.
\newblock {\em Discrete Contin. Dyn. Syst. Ser. B}, 12:905--924.

\bibitem[Zhao et~al., 2010]{Zhao2010}
Zhao, W., Zhang, G., and Ju, L. (2010).
\newblock A stable multistep scheme for solving backward stochastic
  differential equations.
\newblock {\em SIAM J. NUMER. ANAL.}, 48:1369--1394.

\bibitem[Zhao et~al., 2014b]{Zhao2014a}
Zhao, W., Zhang, W., and Ju, L. (2014b).
\newblock A numerical method and its error estimates for the decoupled
  forward-backward stochastic differential equations.
\newblock {\em Commun. Comput. Phys.}, 15:618--646.

\end{thebibliography}
\bibliographystyle{apalike}
\end{document}